\documentclass[1p]{elsarticle}
\pdfoutput=1
\usepackage[utf8]{inputenc}
\usepackage{amsthm}
\usepackage{amsbsy}
\usepackage{amsmath}
\usepackage{graphicx}
\usepackage{amssymb}
\usepackage{lmodern}
\usepackage[pdftex,colorlinks=true,urlcolor=blue,colorlinks=true,urlcolor=blue,citecolor=red,linkcolor=blue]{hyperref}
\usepackage[hyperpageref]{backref}
\usepackage{microtype}

\newcommand{\simInf}[2][+\infty]{\underset{#2 \rightarrow #1} {\sim}}

\newcommand{\convD}{\overset{\mathcal{ D}} {\rightarrow}}
\newcommand{\convP}{\overset{p} {\rightarrow} } 
\newcommand{\Normal}[2]{\mathcal{N}(#1,#2)}

\newcommand{\CharFun}[1] {\mathsf{1}_{#1}}
\newcommand{\InvFun}[1] {#1^{-1}}

\newcommand{\erfc}{\mathrm{erfc}}

\newcommand{\R}{\mathbb{R}}

\newcommand{\Prob}[1]{\mathrm{Prob} \left( #1 \right) }
\newcommand{\Esp}[1]{\mathbb{E} \left[ #1 \right] }
\newcommand{\E}{\mathbb{E}}
\newcommand{\Var}[1]{\mathrm{Var} \left[ #1 \right] }

\newcommand{\qc}{q_c}
\newcommand{\n}{n}
\newcommand{\iid}{\mathrm{i.i.d.}}
\newcommand{\reffig}[1]{\ref{#1}}
\newcommand{\Bias}[1]{\mathcal{B}[#1]} 
\newcommand{\Cov}[2] {\mathrm{Cov} [#1,#2]}

\newcommand{\as}[1]{\overset{\text{a.s}}{#1}}
\newcommand{\Char}[2]{\chi_{#1}\left(#2\right)}
\newcommand{\CharF}[1]{\chi_{#1}}

\newcommand{\logLt}[2][n]{\ln #2 \ll \ln #1}

\graphicspath{
{data/acritical_data/figs/}
{data/acritical_data_n/figs/}
{data/correlation/figs/}
{data/tau_influence/figs/}
{data/linearisation/figs/}
}

\title{Critical moment definition and estimation, \\
 for finite size observation of log-exponential-power law random variables}
\author{Florian Angeletti, Eric Bertin,  Patrice Abry}
\address[ENSL]{Universit\'e de Lyon, Laboratoire de Physique, ENS Lyon, CNRS, UMR 5672
46 All\'ee d'Italie, F-69007 Lyon, France, \\{\tt firstname.lastname@ens-lyon.fr}}

\begin{document}

\bibliographystyle{plain}
\begin{frontmatter}
\begin{abstract}
This contribution aims at studying the behaviour of the classical sample moment estimator, $S(n,q)= \sum_{k=1}^n X_k^{q}/n $, as a function of the number of available samples $n$, in the case where the random variables $X$ are positive, have finite moments at all orders and are naturally of the form $X= \exp Y$ with the tail of $Y$ behaving like $e^{-y^\rho}$. This class of laws encompasses and generalizes the classical example of the log-normal law.
This form is motivated by a number of applications stemming from modern statistical physics or multifractal analysis. 
Borrowing heuristic and analytical results from the analysis of the Random Energy Model in statistical physics, a critical moment $q_c(n)$ is defined as the largest statistical order $q$ up to which the sample mean estimator $S(n,q)$ correctly accounts for the ensemble average $\E X^q$, for a given $n$. 
A practical estimator for the critical moment $q_c(n)$ is then proposed. 
Its statistical performance are studied analytically and illustrated numerically in the case of \emph{i.i.d.} samples. 
A simple modification is proposed to explicitly account for correlation amongst the observed samples. 
Estimation performance are then carefully evaluated by means of Monte-Carlo simulations in the practical case of correlated time series. 
\end{abstract}

\end{frontmatter}



\section{Introduction}
\label{sec:intro}

\begin{figure}
\centerline{ \includegraphics[width=70mm]{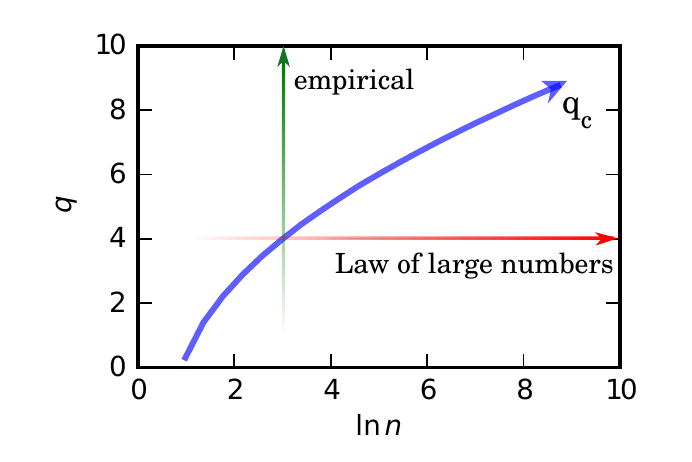} \hspace{2mm} \includegraphics[width=70mm]{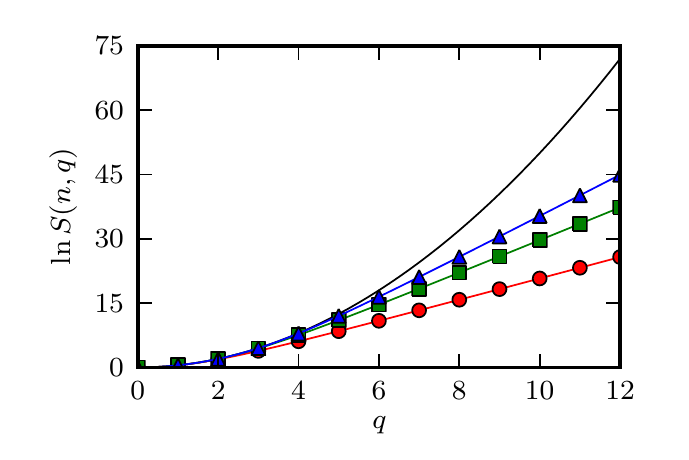}  }
\caption{\label{fig:nqplane} {\bf $(n,q)$ plane and linearization effect.} Left, theoretically, the performance of the sample moment estimator $S(n,q)$ are characterized, for a given $q$, in the limit $n \rightarrow + \infty$, i.e., moving along a horizontal line in the $(n,q)$ plane. Practically, however, the central limit issue consists of determining how large is the order $q$ of the moments that can actually be correctly estimated, from a given observation of finite size $n$, i.e., moving along a vertical line in the $(n,q)$ plane. Right, practically, in the limit $q \rightarrow + \infty$, for a fixed $n$, $\ln S(n,q)$ necessarily behaves as a linear function of $q$, this is illustrated here for a log-normal distribution with three different values of $n$ ( $n=10^2$ (red circle), $n=10^3$ (green square), $n=10^6$ (blue triangle)~; the solid black lines show the theoretical $\ln \E X^q$ ). It shows that both the slope of the linear behavior and the departure from $\ln \E X^q$ depend on $n$. This contribution aims at defining and estimating practically a critical order $q_c(n)$ up to which the sample moment estimator correctly accounts for the ensemble mean $\E X^q$. }
\end{figure}

\noindent {\bf Motivation.} \quad Estimating moments of a given order $q>0$ from a finite size observation of a given times series $\{ X_k, k=1, \ldots, n\}$ appears as both a \emph{natural} and \emph{simple} problem. 
In the general and common case where little a priori is known nor assumed on the data, the classical sample moment estimator for the order $q$,  
$$S(n,q)= \frac{1} {n} \sum_{k=1}^n X_k^{q} $$ 
is natural to use (positive random variables $X$, with finite moments of all orders only are considered, for reasons made explicit below).
Probability theory provides practitioners with valuable results regarding the performance of the sample mean estimator, notably in the asymptotic limit of an infinite observation duration, $n \rightarrow  +\infty $ (cf. e.g., \cite{Billingsley}): 
The weak law of large numbers shows that the estimator is consistent; The central limit theorem further precisely quantifies the asymptotic limit. 
However, such results can also be read as relevant solutions to a rather abstract situation where one is interested in a particular statistical order $q$ and may have at disposal an observation of potentially infinite length. 
In practice, the \emph{natural} situation can often be formulated in a converse way: Practitioners often work hard to obtain an observation of length $n$ and can naturally wonder how large is the order of moments that can actually be estimated correctly given the observed  $\{ X_k, k=1, \ldots, n\}$. 
This dialectic opposition between the theoretical results and the practical questions is sketched in Fig. \ref{fig:nqplane}. 

Even for the case of interest here where all moments are finite, $\E X^q < +\infty$, 
it is well-known that the correct estimation of moments of all orders is impossible, from  any finite size observation.  
Indeed, for a fixed $\n$,
\begin{equation*}
S(n,q) \underset {q \rightarrow +\infty }  {\sim} \frac{1}{n} \max_{k=1,\ldots, n} X_k^{q} \neq \Esp{X^q}.
\end{equation*} 
This behavior that can be referred to as a \emph{linearization} effect, as for fixed $n$, the practical plot, $\log q$ vs. $\log S(n,q) $ systematically appears as a straight line, in the limit $q \rightarrow +\infty$. 
Therefore, a natural question stems: Given an observation of finite size $n$,up to which order moments can actually be correctly estimated~? 
When $X$ belongs to a class of, say, \emph{simple} random variables (defined by the fact that their characteristic function is analytic in $0$), and for \emph{i.i.d.} observations, a bound is given in \cite{KaganNagaev2001}:
  \begin{equation} 
  \label{eq:qcsimple}
  q_c(n) \propto \frac{\ln n}{2 \ln \ln n}.
  \end{equation} 
Though a valuable result, this remains of limited practical use as it is restricted to \emph{i.i.d.} observations belonging to a \emph{simple} class, and, first and foremost, consists of an asymptotical result with no explicit prefactor.

Surprisingly, this question seems to have received little research efforts.
In common real world data analysis, this is however a question of obvious interest.
This is notably the case in scale invariance analysis, increasingly used in many applications of very different natures (cf.~e.g., \cite{Abry2002} and references therein for a review), where scaling is measured via the dependence with respect to the analysis scale $a$ of the moments of order $q$ of multiresolution quantities computed from the observation $\{ X_k, k=1, \ldots, n\}$.
Such an issue also appears, though in a less obvious way, in the physics of disordered systems,
a representative example of which being the so-called Random Energy Model (REM) \cite{Derrida1981,mezard1984}, where the partition function
\begin{equation}
Z(n,\beta) = \sum_k \exp (\beta Y_k),
\end{equation}
(a fundamental quantity in the statistical physics framework)
needs to be evaluated, with $\beta = 1/k_B T$, $T$ the temperature of the system and $k_B$ the Boltzmann constant. The variables $Y_k$ are (the opposite of) random energies, drawn from
a Gaussian distribution in the original version of the model \cite{Derrida1981}.
Such a choice leads to a log-normal statistics for the variables $X_k=\exp(\beta Y_k)$.
The REM has however been extended to other random variables, such as log-exponential-power laws, obtained by choosing a distribution of $Y_k$ with a tail $p(y) \sim exp(-c y^\delta)$, with
$\delta>1$ and $c>0$, when $y \to +\infty$ \cite{BouchMez97}.
This similarity between the REM and the moment estimator issue raised above stems from
the close resemblance between the partition function $Z(n,\beta) $ and the sample moment estimator $S(n,q)$ for random variables $X_k = \exp Y_k$, with $Y_k$ belonging to the \emph{simple} class.

Another interesting example stems from multifractal analysis. Indeed, popular multifractal processes are defined from the celebrated Mandelbrot multiplicative cascades (cf. e.g., \cite{Mandelbrot1974,frisch1995}), where the multipliers entering the construction of the cascade are also chosen as $X = \exp (-Y)$. \\

\noindent {\bf Goals, Contributions and Outline.} \quad In this context, the present contribution first aims at defining a critical order $q_c(n)$ up to which moments of order $q$ can be estimated from a given finite size observation of length $n$, belonging to the class of positive RV $X = \exp Y$, where the class of RV to which $Y$ belongs is precisely defined in Section \ref{sec:RV}.
The definition of the critical order $q_c(n)$ is based on a detailed heuristic argumentation borrowed from the analysis of the REM \cite{Derrida1981,mezard1984} and is also strongly grounded on the thorough theoretical asymptotic analysis conducted in \cite{BenArous2005}. 
This construction is detailed in Section \ref{sec:analyse}. 
Second, a procedure for the practical estimation of $q_c(n)$ from a finite duration observation is devised and its performances are studied theoretically in the case of \emph{i.i.d.} observations (cf. Section \ref{sec:estim}). 
Third, the estimation procedure is slightly modified to account for potentially dependent observations and its performances are thoroughly studied by means of Monte-Carlo simulations  (cf. Section \ref{sec:cov}). 
The {\sc Matlab} practical procedure implementing estimation will  be made available publicly at the time of publication. 

\section{Class of random variables}
\label{sec:RV}

\noindent {\bf Log-normal.} \quad As mentioned in Section \ref{sec:intro} above, the typical example of interest here is that of log-normal random variables $X$. 
In that case, $Y = \ln X $ simply consists of a normal random variable, whose cumulative distribution $F_Y$ satisfies $\ln (1 -F_Y(y)) \approx y^2/2$ for $y\rightarrow< +\infty$. 
In the present contribution, following extensions naturally encountered in both the REM and Multifractal contexts and the theoretical framework envisaged in \cite{BenArous2005}, a much larger class of random variables is considered by extending the log-normal case in two ways: The power law behavior $y^2$ is extended to $y^{\rho}$. The power-law behaviour is defined through asymptotic properties which capture the essential characteristics of exact power-law.\\

\noindent {\bf log-exponential-power law.} \quad The class of \emph{log-exponential-power law} random variables $X$ is defined as follows: \begin{equation}
X = \exp Y,
\end{equation}
where the cumulative distribution function of $Y$, $F_Y$ is assumed to be a continuous function that satisfies: 
\begin{equation} 
\label{eq:h}
 1-F_Y(y) = e^{ -h(y)}. 
\end{equation} 
Furthermore, the function $h$ is required to verify:
\begin{equation}
\label{eq:hform}
h(y)=L(y) y^\rho, \, \makebox{ with } \rho > 1, 
\end{equation}
where $L$ is a slowly varying function (cf. e.g., \cite{Beran1994}), i.e.,  a function satisfying $$ \forall t, \lim_{y \rightarrow +\infty} {L(ty) \over L(y)} = 1. $$
In addition, for technical reasons made clearer in Section \ref{sec:analyse}, this class is slightly restricted to twice-differentiable functions $L$ satisfying: 
\begin{eqnarray}
\label{eq:hconda}
\lim_{y \rightarrow +\infty} y {L'(y) \over L(y)}  & = & 0,\\
\label{eq:hcondb}
 \lim_{y \rightarrow +\infty}  \left(y {L'(y) \over L(y)}\right)'  & = & 0.
\end{eqnarray}
These condition impose a certain form of regularity on the derivative of $L$. Notably Eq.(\ref{eq:hconda}) excludes slowly varying function with increasingly fast oscillations. Similarly Eq.(\ref{eq:hcondb}) constraints the behaviour of $L''$ with respect to $L$ and $L'$. 

\noindent {\bf Consequences.} \quad From the definition above, a number of properties for the random variables $Y$ and hence $X$ can be derived. 
The probability density function $p_Y$ of $Y$ reads: 
\begin{equation} 
p_Y(y) = h'(y) e^{-h(y)}. 
\end{equation}
All moments of X are finite and can be written as:
\begin{equation}
\label{eq:EXq}
\Esp {X^q} = \int_{-\infty}^{+\infty} h'(y) e^{q y-h(y)} dy = \int_{-\infty}^{+\infty} h'(y) e^{q y-y^\rho L(y)} dy.
\end{equation}
Nevertheless, the moment generating function of $X$,
\begin{equation}
\Esp {e^{tX}} = \int_{0}^{+\infty} h'(y) e^{t\exp(y)-L(y)y^\rho} dy,
\end{equation}
 is finite on $\R^-$ only and its characteristic function is hence not analytical in 0 (see \cite{Lukacs} for more detailed consequences).

\noindent {\bf Interpretation.} \quad Eq.(\ref{eq:hform}) together with Conditions (\ref{eq:hconda}) and (\ref{eq:hcondb}) essentially control the behavior of the right tail of the distribution of $Y$ and hence define a broad class of random variables $X$, whose moments are all finite, but with characteristic function not analytical at zero. 
The interest of that class of random variables hence stems from its lying between the \emph{simple} class of random variables, whose moments are all finite and whose tail decreases exponentially fast (whose 
 characteristic function is hence analytic at zero, e.g., the absolute value of normal random variables) and the so-called \emph{heavy tailed} class, whose moments remain finite only up to a given order (such as e.g., the Pareto distribution).
The parameter $\rho$ plays a key role in studying that evolution from \emph{simple} to \emph{heavy tail}: 
while, strictly speaking, both classes are excluded from the class studied here, they can be approached by letting $\rho \rightarrow +\infty $ and $\rho \rightarrow 1$, respectively. \\

\noindent {\bf Local power law exponent.} \quad In Definition (\ref{eq:hform}), the parameter $\rho$ is naturally referred to as the power-law exponent. 
However, in the course of the study reported below, a \emph{local} power-law exponent, defined as follows, will play a crucial role, where \emph{local} refers to a pre-asymptotic behavior, i.e., for finite $y$: 
\begin{equation} 
\label{eq:rhol}
\rho_l(y) = \frac{\partial \ln h(y)}{\partial \ln y}=  \frac{y h'(y)}{h(y)}.
\end{equation}
It is straightforward to verify that Conditions (\ref{eq:hconda}) and (\ref{eq:hcondb}) imply:
\begin{eqnarray} 
\lim_{y\rightarrow +\infty} \rho_l(y) & = & \rho, \\
\lim_{y\rightarrow +\infty }\rho_l'(y) & = & 0.
\end{eqnarray}

\noindent {\bf Examples.} \quad In the numerical analyses conducted in Sections \ref{sec:estim} and \ref{sec:cov} below, two explicit examples of $\rho$-parametrized families of such random variables $X$ are used. 
\begin{itemize}
\item The log-Weibull distribution defined as:
\begin{equation} 
F(y) = 1- e^{-y^\rho}, \, \makebox{ with } \rho > 1,y>0.  
\label{eq:def:lW}
\end{equation}
This notably implies that the local power exponent is constant 
\[  \forall y, \rho_l(y) \equiv \rho. \]
\item The strict log-exponential-power law distribution defined as: 
\begin{equation} 
\label{}
F_Y(y)= \frac{1}{2 \Gamma(1+1/\rho)}\int_{-\infty}^{y} e^{-|t|^\rho} dt, \, \makebox{ with } \rho > 1
 \end{equation} 
where $\Gamma$  denotes the gamma function \cite{Royden}. 
This implies that its probability density function reads: 
\begin{equation} 
p_Y(y) = \frac{\rho}{2 \Gamma(1+\frac{1}{\rho}) } e^{-|y|^\rho}.
\label{eq:def:lep} 
\end{equation}
It can be shown that using the asymptotic properties of the incomplete $\gamma$ functions that
\begin{equation*} \lim_{y-\rightarrow +\infty} \rho_l(y)=\rho. \end{equation*}
 
\item  The log-normal distribution corresponds to a special case of the strict log-exponential-power law family, with $\rho=2$. 
The following relations can be derived: 
\begin{eqnarray*}
h(y) & =  & \ln 2 - \ln \erfc (\frac{y}{\sqrt 2}), \\
\rho_l(y)& =  & 2 \frac{y \exp(-y^2/2)} {\sqrt{2 \pi}\erfc (y/\sqrt 2)}.
\end{eqnarray*}
The asymptotic form of the complementary error function $\erfc$, in the limit $y\rightarrow + \infty$ : $\erfc(y) \sim \frac{\exp(-y^2)}{\sqrt{2\pi} y } $ 
yields: 
\begin{eqnarray*}
h(y) & \sim  & \ln 2 + \frac{y^2}{2} + \ln (\sqrt{2\pi} y),  \\
\lim_{y\rightarrow + \infty} \rho_l(y) & =  & 2 .
\end{eqnarray*}
\end{itemize}

\section{Critical order: combining dominant contribution to finite size effects and truncated moments.}
\label{sec:analyse}

This section aims at unveiling a change in the behaviour of the sample mean moment estimator,
\begin{equation}
S(n,q)= \frac{1} {n} \sum_{k=1}^n X_k^{q},
\end{equation} 
in the $(q,n) $ plane,  for a critical curve  $q_c(n)$, that hence depends on the observation sample size $n$. 
This is obtained by combining a a moment dominant contribution argument to a truncated moment finite size effect analysis. 
This analysis is strongly inspired by those commonly conducted in statistical physics, notably for the study of the Random Energy Model \cite{Derrida1981} and 
closely resembles a glass transition. 

\subsection{Moment dominant contribution}
\label{sec:mdc}

The behavior of the moment of order $q$ of $X$ can first be related to the tail of its distribution and hence to $h(y)$. 
Eq.(\ref{eq:EXq}) can be rewritten as: 
\begin{equation} 
\label{eq:EXqbis}
\Esp{X^q}= \int_0^{+\infty} e^{qy -h(y) + \ln h'(y)} dy. 
\end{equation}
Using a saddle-point technique \cite{DeBruijnAsymptotic}, in the limit $q \rightarrow +\infty$, enables one to evaluate the integral above as (technical proof detailed in \ref{app:SaddlePoint}): 
\begin{equation} 
\label{eq:momsaddlepoint}
\Esp{X^q} \underset{q \rightarrow +\infty}{\sim} 
 \sqrt{\frac{2 \pi} {h''(y^*)+ (\ln h')''(y^*)} }
 \exp \left( qy^* -h(y^*) + \ln h'(y^*) \right),
\end{equation}
where  $y^*$ is determined by
\begin{equation}
\label{ystarq}
 q - h'(y^*) + \frac{h''}{h'} (y^*)=0. 
 \end{equation} 

Qualitatively, $y^*$ denotes the value of $y$ 
that essentially contributes to the moment of order $q$, $\Esp{X^q}$. 
 Eq. (\ref{ystarq}) implicitly defines a relation $y^*(q)$, that depends only on the precise details of the distribution of $Y$. 
Because Condition (\ref{eq:hcondb}) imposed in Section \ref{sec:RV} precisely aims at constraining its amplitude, the term $h''/ h'$ can be read as a correction, depending of the fine local behavior of $p_Y$. 
Indeed, under this condition, one can show that: 

\begin{equation*}
\lim_{y\rightarrow +\infty} \frac{h''(y)}{h'(y)}=0,
\end{equation*}
which hence implies that:
\begin{equation}
h'(y^*) \approx q.
\end{equation}

\subsection{Finite size effect: truncated moments and sample moment estimator}
\label{sec:fs}

Let us consider a finite size observation of $ n $ \emph{i.i.d.} observations  $\left\{Y_k, k \in 1 \dots \n \right\}$. 
Its samples lie between empirical minimum and maximum $m_k=\min \{ Y_k\}$ and $M_k=\max \{ Y_k \}$. 
A finite size observation thus explores only a bounded range of values in the existence domain of $y$. 
When $n \rightarrow +\infty$, the bounds of this range naturally extend to cover the entire domain of definition of $F_Y$.
However, for finite $n$, this bounded range is likely to affect the properties of the sample mean estimator for a given $q$, if the information needed to evaluate this moment (qualitatively $y^*$)  lies beyond the observed $M_k$. 
This section aims at detailing and quantifying this statement as well as at drawing its consequence in terms of sample moment estimator properties. 

\subsubsection{Finite size accessible range} 

Let us first define a frontier beyond which the probability to observe a sample $Y_k=y$ within a  finite size observation is low. 
Let $\tau $ denote an arbitrary positive constant, that controls the probability that $ \makebox{max}\{ Y_k, k=1, \ldots, n\} $ exceeds a frontier $y^\dagger_\tau$: 
\begin{equation*} 
\Prob{\forall k=1, \ldots, n, Y_k < y^\dagger_\tau} = e^{-\tau}. 
\end{equation*}
Using the independence of the $\left\{Y_k, k \in 1 \dots \n \right\}$ and the relation above yields, $1-F(y^\dagger_\tau) = \frac{\tau}{n} + o (\frac{1}{n}) $, or equivalently, $ h(y^\dagger_\tau) = \ln \tau + \ln n+o(1) $. 
For the class of random variables considered here (cf. Eq.(\ref{eq:hform})), for $n \rightarrow +\infty$, one has: $ y^\dagger_\tau \sim \left( 1 + \frac{\ln \tau}{\ln n} \right)^{\frac{1}{\rho}} y^\dagger_1 $. 
This indicates that it is natural to simply consider $y^\dagger(n)$, defined as:
\begin{equation}
\label{eq:ydag}
 h(y^\dagger(n))= \ln n. 
\end{equation} 
which quantifies the upper bound of the tail of $y$ actually likely to be observed from an observation of finite size $n$. 
Let us note that $y^\dagger(n)$ depends both on $n$ and on the details of the tail of the distribution of $Y$.

\subsubsection{Truncated moments} 

Essentially, the finite size analysis above states that, for a given observation of finite size $n$, the domain of existence of $Y$, $\R$, must be split into two distinct subranges: 
\begin{itemize}
\item[-] { $(-\infty, y^\dagger)$, a reachable range which contain most observed samples,\footnote{For positive moments, the positives value of $Y$ play a dominant role, so one can focus of the positive frontier $y^\dagger$. However, a natural extension to negative moments would introduce the negative counterpart $y^\dagger_-$ of $y^\dagger$ which would play an exactly symmetric role.}}
\item[-] $[y^\dagger, +\infty)$, an unreachable range, within which the probability of observed samples is very low. 
\end{itemize}

This leads us to define the so-called truncated moments: 
\begin{equation} 
M_T(n,q)= \int_0^{y^\dagger(n)} e^{qy -h(y) + \ln h'(y)} dy. 
\end{equation}
This truncated moment provides an heuristic estimate of the typical value of the moments, where the contributions of atypical term is removed. 
\subsection{Sample moment estimator and truncated moment} 

The sample moment estimator $S(n,q) = (1/n) \sum_{i=1}^{n} X_i^q $ for an observation of $n$ $\iid$ random variables $X_i = e^{Y_i}$, in the class of power law exponential, of interest here, can be related to the truncated moments as follows. 
Let $q(n) $ denote a function of $n$ that satisfies: 
\begin{equation}
\label{eq:qn}
\exists \epsilon>0, q(n) \ll (\ln n)^{2-1/\rho-\epsilon},
\end{equation}  
then, as $ n \rightarrow +\infty$,
\begin{equation} 
\label{eq:SMT}
\frac {\ln S(n,q(n))} {\ln n} \overset {\text{a.s}} {\rightarrow} \lim_{n \rightarrow +\infty} \frac{\ln M_T(n,q(n))}{\ln n}.
\end{equation} 
Essentially, this means that, asymptotically ($ n \rightarrow +\infty$), on condition that $q$ does not grow too fast (cf. Condition (\ref{eq:qn})), the (ln of the) sample moment estimator converges to the (ln of the) truncated moment. 
The complete proof of this results is detailed in \ref{app:LogConv}. 
In essence,  it relies on the specific properties of the class of random variables studied here and closely exploits in spirit the analogy with the REM, extending reasonings such as those in \cite{3}. 


\subsection{Phase transition and critical moment}

\subsubsection{Critical frontier}

Combining the moment dominant contribution evaluation (cf. Section \ref{sec:mdc}) to finite size truncation effects (cf. Section \ref{sec:fs}) yields two distinct situations depending on the relative positions of $y^*(q) $ and $  y^\dagger(n)$, and hence separated by:
\begin{equation} 
\label{eq:ydagystar}
y^\dagger(n)=y^*(\qc).
\end{equation}
This equality implicitly defines a frontier that splits the $(q,n)$ plane, introduced in Section~\ref{sec:intro}, into two zones corresponding to different behaviours of the truncated moments (in the REM framework, this is a phase transition). 
This frontier can either be formulated as a critical order $q_c(n)$ or as a critical size $n_c(q)$. \\

For $y^*(q)  \leq  y^\dagger(n)$ (equivalently for $q \leq q_c(n)$), the saddle-point argument enables, as detailed in \ref{app:SaddlePoint}, the evaluation of $M_T$ as 
\begin{equation} 
\ln M_T(n,q) =   q y^*(q) -h(y^*(q)) + \ln h'(y^*(q)),  
\end{equation}
which by comparison to Eq.~(\ref{eq:momsaddlepoint}), is hence shown to correspond to $\ln \Esp{X^q}$. 
By nature, the saddle-point argument implies $q \rightarrow + \infty$, which combined to $q \leq q_c(n)$ yields that the result above is valid for $q \rightarrow + \infty$, $n \rightarrow + \infty$, with $q \leq q_c(n)$. 

Conversely, when  $y^*(q) \geq y^\dagger(n)$  (equivalently for $q \geq q_c(n)$), the dominant contribution to $M_T(n,q)$ is no longer controlled by $y^*(q)$ but instead by the border of the integration domain $y^\dagger(n)$. 
Therefore, the saddle-point argument now yields, in the limits $q \rightarrow + \infty$, $n \rightarrow + \infty$, with 
$q \geq q_c(n)$ : 
\begin{equation} 
\ln M_T(n,q) =    q y^\dagger(n) -h(y^\dagger(n)) + \ln h'(y^\dagger(n)). 
\end{equation}
This clearly shows that $\ln M_T$ depends linearly on q, for large $q$. 


\subsubsection{Sample moment estimator}

The almost sure convergence established in Eq.~(\ref{eq:SMT}), and proven to be valid regardless of the relative values of $y^*(q) $ and $  y^\dagger(n) $, indicates that the phase transition is also transferred to the behaviour of the sample moment estimator, which reads, 
in the limits $q \rightarrow + \infty$, $n \rightarrow + \infty $, with  $q_M(n)=(\ln n)^{2-1/\rho-\epsilon}$ satisfying Condition~(\ref{eq:qn}): 
\begin{equation}
\ln  S(\n, q ) \approx
\begin{cases}
 \ln \Esp{X^q}  & q  \leq q_c(n)\\
(\ln n \rho) \frac{q}{q_c}  & q_c(n) \leq q \leq q_M(n).\\
 \end{cases} 
 \end{equation}
As can be seen comparing Eqs. (\ref{eq:qn}) and (\ref{eq:qc}) below, for $\epsilon \in (0,1)$  $q_c(n) \ll q_M(n)$. 
Therefore, for a large range of $q$s, the sample moment estimator $S(n,q)$ correctly estimates the ensemble average $\E X^q$ for $q \leq q_c(n)$, while $ \ln S(n,q)$ behaves linearly in $q$ when $q \geq q_c(n)$. 

\begin{figure} [ht]
\centerline{\includegraphics[width=70mm]{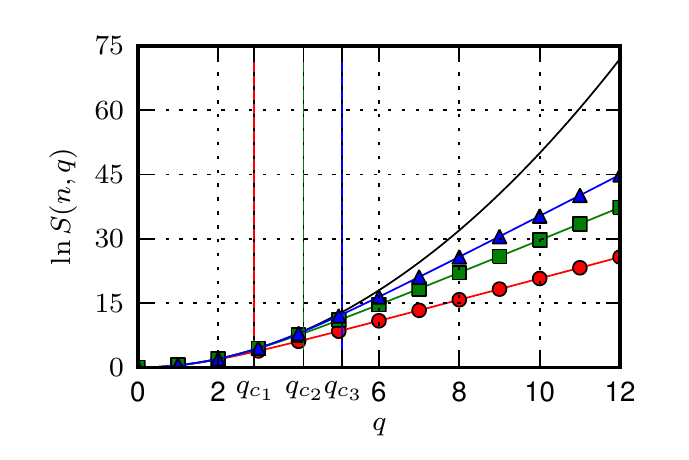} \hspace{2mm} \includegraphics[width=70mm]{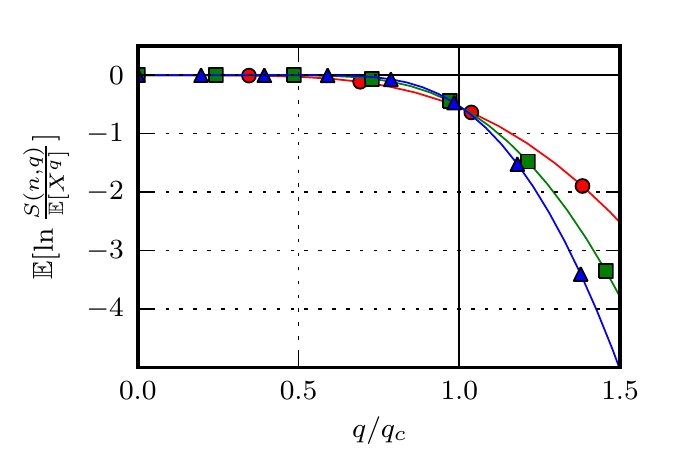} }
\caption{\label{fig:qc:convergence}{\bf Quantitative predictions for the linearization effect.} Left, $ \E \ln S(n,q)$ as a function of $q$, with predicted values of $q_c$, for three different values of $n$: $n=10^2$ (red circle), $n=10^3$ (green square), $n=10^6$ (blue triangle)~; the solid black lines show the theoretical $\ln \E X^q$. Right, ${\Esp{\ln S(n,q)}} / \Esp{X^q} $ as a function of ${q / \qc}$, for the same $3$ different values of $n$.  $ \E \ln S(n,q)$ is obtained from Monte-Carlo simulations, with a Log-Normal distribution. 
These plots clearly validate the relevance of the defined critical moments $q_c(n)$.}
\end{figure}

This phase transition is illustrated in Fig.~\ref{fig:qc:convergence} from numerical simulations.
A large number of independent realizations of $n$ \emph{i.i.d.} log-normal RV $X$ enable to compute an average value of  $\ln S(n,q)$, for various $n$ and $q$.
This average can further be compared to $\ln \Esp{X^q}$.
This clearly validates that for $q< q_c(n)$,  $\ln S(n,q)$ corresponds to $\ln \Esp{X^q}$, while $\ln S(n,q)$ becomes significantly lower than $\ln \Esp{X^q}$ for $q> q_c(n)$ and behaves linearly in $q$. 
Results obtained for various $n$ can be superimposed by plotting $\ln S(n,q)$ and $\ln \Esp{X^q}$ as a function of $q/q_c(n)$ rather than as a function of $q$, where $q_c(n)$ is obtained by solving numerically Eq.~	(\ref{eq:ydagystar}) for each $n$. 
  
\subsubsection{Critical moment}

Combining Eq.~(\ref{eq:ydagystar}) with Eqs.~(\ref{ystarq}) and (\ref{eq:ydag}) yields:  
\begin{equation} 
\label{eq:ydagystarbis}
\frac{\qc}{\ln n} = \frac{1}{y^\dagger} \frac{y^\dagger h'(y^\dagger)}{h(y^\dagger)}  - \frac{h''(y^\dagger)}{h(y^\dagger) h'(y^\dagger)  }
\end{equation}
Given Conditions (\ref{eq:hconda}) and (\ref{eq:hcondb}), the second term on the right side of the equation above can be shown to be negligible compared to the first one. 
Moreover,  on the right-hand-side of Eq. (\ref{eq:ydagystarbis}), the local tail exponent $\rho_l$ can be recognized: 
\begin{equation}
\label{eq:rhoe}
 \frac{y^\dagger h'(y^\dagger)}{h(y^\dagger)}  = \rho_l (y^\dagger(n)). 
\end{equation}

Further defining,
\begin{equation}
\label{eq:theta}
\theta(n) = \frac{\ln n}{y^\dagger(n)},
\end{equation}
 $q_c(n)$ can be approximated as:  
\begin{equation}
\label{eq:qc:1}
q_c(n) \simeq  \rho_l(y^\dagger(n)) \theta(n),
\end{equation}
which will prove of particular interest for the practical estimation of $q_c(n)$. 

As an interesting corollary, the growth rate of $q_c$ as a function of $n$ and $\rho$ can be evaluated.
Eq. (\ref{eq:hform}) implies that there exists a slowly varying function ${\cal L}$, closely related to $L$, cf. \cite{BinghamRegularVariation},  such that:  
\begin{equation*} 
{\InvFun{h}(\ln n) = {\cal L}(\ln n) (\ln n)^{1/\rho} }, 
\end{equation*}
Substituting this into Eq. (\ref{eq:qc:1}) leads to 
\begin{equation} 
\label{eq:qc}
\qc(\n) = \rho {\cal L}(\ln n) (\ln n)^{1-\frac{1}{\rho}}.  
\end{equation}
When $ \rho \rightarrow + \infty $, this is reminiscent of Eq. (\ref{eq:qcsimple}), recalled in Section \ref{sec:intro} and valid for the \emph{simple} class of random variables. 
However, when $\rho$ decreases and approaches $1$, this critical moment $q_c$ increases far more slowly with $\n$. 
As an example, for the log-normal case, $ \rho = 2$ implies that $\qc(n) \propto \sqrt{\ln n} $. 
 
\section{Estimation}
\label{sec:estim}

This section aims at defining an estimation procedure for $q_c(n)$ that can actually and practically be used from a finite sample size observation $\{ X_1, \ldots X_n \}$. 
Also, its statistical performance are studied both theoretically and numerically. 

In Section \ref{sec:analyse} above, it has been shown (cf. Eq.~(\ref{eq:qc:1})) that $q_c(n) $ can naturally be split into two components $\theta(n)$ and $\rho_l(y^\dagger(n))$. 
An estimation procedure is hence devised and studied independently for each components in Sections~\ref{sec:thetaest} and \ref{sec:rhoest}, respectively, yielding an estimate of $q_c$ as a product of the estimates (cf. Section~\ref{sec:qcest}). 

\subsection{Estimation of $\theta$}
\label{sec:thetaest}

\subsubsection{Definition}

A straightforward calculation not reported here shows that $ y^\dagger(n)= \InvFun h(\ln n) = \InvFun F_Y(1 - \frac{1}{n})  $, which indicates that 
the definition of $\theta$ (cf. Eq. \ref{eq:theta}) can be rewritten in terms of the $(1-1/n)$-th quantile of the distribution of $Y$:
\begin{equation} 
\label{eq:deftheta}
\theta = \frac {\ln n} {\InvFun F_Y(1 - \frac{1}{n}) } .
\end{equation}
The estimation of  $\theta$  can, therefore, essentially be recast into that of (the inverse of) the  $(1-\frac{1}{n})$-th quantile of $Y$.

Classically,  the $p$-th quantile is estimated, from an observation of sample size $n$, using the order statistics of rank  $\left[ p . n \right]$ \cite{Billingsley}. 
In the present context, this amounts to using the largest value of the observation, 
\begin{equation} 
M_{Y,n} = \makebox{ Max }\{ Y_1, \ldots, Y_n \}, 
\end{equation}  
as an estimate of the $(1-\frac{1}{n})$-th quantile.
For the class of random variables $Y$ considered here, $M_{Y,n}$ belongs to the domain of attraction of Gumbel law \cite{EmbrechtsExt,GumbelExt,GalambosExt}.
This implies that there exists  a sequence $a_n \geq 0$ such that $[M_{Y,n} - \InvFun F_Y(1-\frac{1}{n})]/{a_n}$ converge in distribution towards a Gumbel distribution with a cumulative function $F(x)=exp( -e^-x)$.
Consequently,
\begin{equation} 
\lim_{n\rightarrow +\infty} \Esp {\frac{M_{Y,n} - \InvFun F_Y \left(1-\frac{1}{n} \right)}{a_n}}= \gamma, 
\end{equation}  
where $\gamma $ is the Euler-Mascheroni constant \cite{Royden}.
This straightforwardly indicates that $M_{Y,n}$ constitutes a biased estimator for $ \InvFun F_Y(1 - \frac{1}{n}) $ for finite sample size and that it is not, in general, asymptotically consistent. 

To overcome this drawback, the estimation procedure can be refined by involving not only $M_{Y,n}$, but the $k$-largest values in the observation $\{ Y_1, \ldots, Y_n \}$. 
Therefore, let $\{ Y_{i,n}, i=1, \ldots, n \} $ denote the ordered list (in descending order) of the observation $\{ Y_{k}, k=1, \ldots, n \} $ and let $ k_\theta(n)$ be a function of $n$ satisfying (for reasons made clear later):
\begin{equation} 
\label{eq:kn}
\lim_{n \rightarrow + \infty} \frac{k_\theta(n)}{n}=0. 
\end{equation} 
A parametrized collection of estimators of $\theta$ is defined as (with $k \equiv k_\theta(n)$): 
\begin{equation} 
\label{eq:thetahat}
\hat{\theta }^{(k)} = \frac{\ln n}{\Omega_k},
\end{equation}
where ${\Omega_k}$ consist of linear combinations of the order statistics: 
\begin{equation} 
\label{eq:omega:def} 
\Omega_k= \sum_{i=1}^k \alpha_i Y_{i,n}, \makebox{ with } \sum_{i=1}^k \alpha_i \equiv 1.
\end{equation}  
Obviously, the case $k = 1$ amounts to using $M_{Y,n}$ only. 

\subsubsection{Performance: theoretical analysis}

To study the performance of $\hat{\theta }^{(k)}$, it is needed to analyze the statistical properties of ${\Omega_k}$, which can be rewritten as 
\begin{equation*} 
\label{eq:omegarewrite}
\Omega_k =   \InvFun F \left(1 - \frac{1}{n} \right) + a_n \sum \alpha_i U_i, 
\end{equation*}
where
\begin{equation*} 
U_i = \frac{Y_{i,n} -   \InvFun F(1 - \frac{1}{n})}{a_n}. 
\end{equation*}
In the limit $n \rightarrow +\infty $, the $\{ U_i, i=1, \ldots, k=k_\theta(n) \} $ converge in distribution ($\overset{\mathcal{D}}{\rightarrow}$) toward a random vector:
\begin{equation*}
\left(U_1, \dots, U_k \right)\overset{\mathcal{D}}{\rightarrow}\Lambda^{(k)} \equiv  \left(G_1, \dots G_k \right),
\end{equation*}
whose joint distribution reads \cite{EmbrechtsExt,GalambosExt}: 
\begin{equation} 
p_{\Lambda_k} (g_1, \dots g_k)=\CharFun{g_{1}>\dots>g_{k}}\exp\left(-e^{-g_{k}}-\sum_{i}g_{i}\right).
\end{equation}  
The change of variable $\Lambda^{(k)} \equiv \{G_1, \dots G_k \} \rightarrow \Lambda'^{(k)}\equiv   \{ \{\Delta_i =i(G_{i}-G_{i+1)}, i=1, \ldots, k-1 \}, G_k \}$ yields: 
\begin{equation} 
p_{\Lambda'_k} (\delta_1, \dots \delta_{k-1},g_k )=\exp\left(-e^{-g_{k}}-kg_k-\sum_{i<k}\delta_{i}\right), 
\end{equation} 
which shows that asymptotically  ($n \rightarrow + \infty$), the $ \{ \Delta_i, i=1,\ldots,k-1 \}$ are independent.\\

From these definitions, a number of properties of ${\Omega_k}$, for $k \geq 2$, have been derived and are presented below. 

\newtheorem{theor_Omega}{Proposition}
 \begin{theor_Omega}
 \label{theor_Omega}
On condition that it is applied to a random vector that exactly follows $p_{\Lambda_k}$, 
 the estimator $\Omega_k$ that has minimal variance under the constraint that it is unbiased for finite $n$ (i.e., $\E \Omega_k = \InvFun F(1-\frac{1}{n})$ 
 or, equivalently,  $\E \sum_{i=1}^{k} \alpha_i G_i  \equiv 0$ ) is obtained by setting $\{ \alpha_i, i = 1 \ldots, k \} $ to: 
\begin{equation} 
\label{eq:alphadef}
\begin{cases}
\forall i, i\neq k & \alpha_i= \frac{\gamma- \sum_{l=1}^k \frac{1}{l}}{k-1} \\
 & \alpha_k=1- (k-1)\alpha_1 \\
\end{cases}
 \end{equation}
 \end{theor_Omega}
 
This stems directly from the following results (whose proofs are postponed to \ref{sec:theta_appendixa})
\begin{eqnarray} 
\Esp{G_k} &=& \gamma - \sum_{j=1}^{k-1} \frac{1}{j}, \\
\Var{G_{k}} &=& \frac{\pi^{2}}{6}-\sum_{j=1}^{k-1} \frac{1}{j^2}, \\
\Esp { \Delta_j }&=& 1, \\
 \Var{\Delta_j} &=& 1. 
\end{eqnarray} 

 \begin{theor_Omega}
When $n \rightarrow + \infty$, $k \rightarrow +\infty$ with $ \frac {k} {n} \rightarrow 0 $, $\Omega_k$ is consistent.
 \end{theor_Omega}
 
This is a direct consequence of the combination of  Proposition \ref{theor_Omega} above with the fact that for all random variables in the class studied in this contribution, the random vector $\{ U_i, i=1, \ldots, k \} $ converges in distribution towards $ \Lambda^{(k)} $. Proof is further detailed in \ref{sec:theta_appendixb}. \\

\begin{theor_Omega}
 When $n \rightarrow + \infty$, $k \rightarrow +\infty$ with $ \frac {k} {n} \rightarrow 0 $, $\Omega_k$ is asymptotically normal.
 \end{theor_Omega}
When applied to a random vector that exactly follows $p_{\Lambda_k}$, $\Omega_k$ can be split into
 \begin{equation*}
 \Omega_k = G_k + \alpha_k \sum_i \Delta_i.
 \end{equation*}
Using the fact that the vector $ \{ \Delta_i, i=1,\ldots,k-1 \} $ consists of unit variance \emph{i.i.d.} random variables yields
\begin{equation*}
 \frac{\Omega_k+\Esp {G_k}-G_k} {(\zeta(1;k-1)-\gamma)^2} \convD \Normal{0}{1} \text{ where } \zeta(s;n)=\sum_{k=1}^n \frac{1}{k^s}. 
 \end{equation*}
 and
 \begin{equation*} 
 \Esp{G_k}-G_k \convP 0.
 \end{equation*} 
This implies: 
\begin{equation*} 
\frac{\Omega_k} {(\zeta(1;k-1)-\gamma)^2} \convD \Normal{0}{1}.
\end{equation*}
This shows that $\Omega_k$ is normal for finite $n$ when applied to a random vector that exactly follows $p_{\Lambda_k}$. 
Asymptotic normality in general follows from the convergence in distribution of $\{ U_i, i=1, \ldots, k=k_\theta(n) \} $ towards $ \Lambda^{(k)} $. 
Proof is further detailed in \ref{sec:theta_appendixb}.\\

The consistency of $\theta^{(k)}$ immediately follows from that of $\Omega_k$. 
However, it should be noted that $\theta^{(k)}$ is not asymptotically normal.  

These analytical results are obtained with choices $k_\theta(n)$ satisfying the constraint in Eq.~(\ref{eq:kn}). 
Choosing precisely $k_\theta(n)$ is however an intricate question that can be formulated into a classical bias-variance trade-off one: 
Large $k$ should yield low variances while small $k$ are needed to ensure convergence toward a law of extremes ($\Lambda^{(k)}$) and hence a low bias. 
This trade-off is now studied by means of Monte-Carlo simulations.

\subsubsection{Performance: numerical analysis}
\label{sec:thetanum}

 \begin{figure}
\centerline{}
\centerline{\includegraphics[width=60mm]{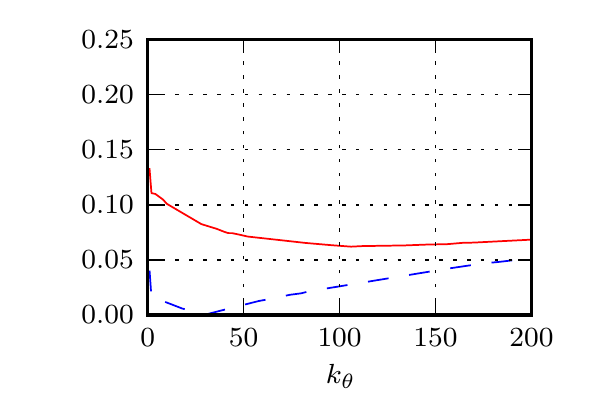} \hspace{1mm} \includegraphics[width=60mm]{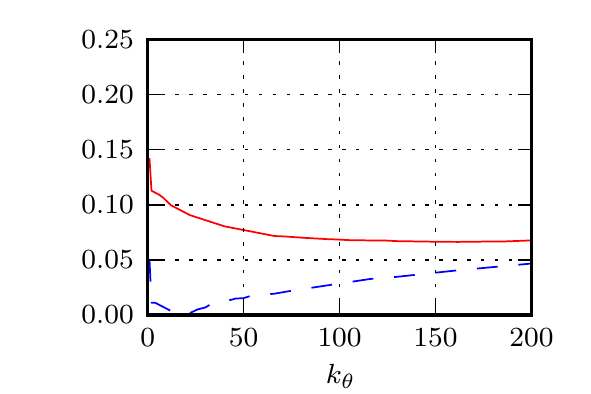}} 
\centerline{ \includegraphics[width=60mm]{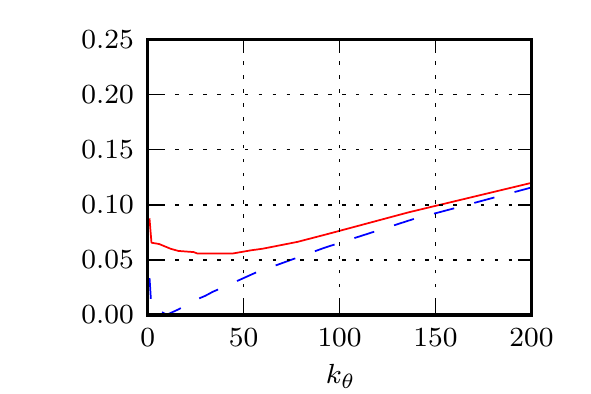} \hspace{1mm} \includegraphics[width=60mm]{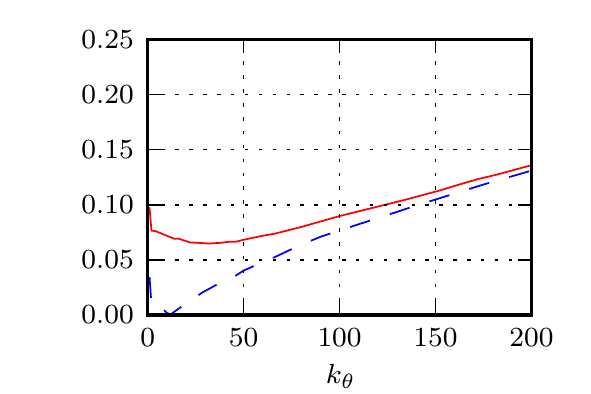} } 
\centerline{ \includegraphics[width=60mm]{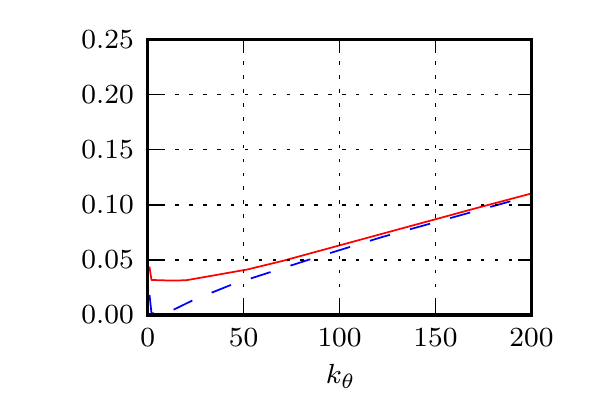} \hspace{1mm} \includegraphics[width=60mm]{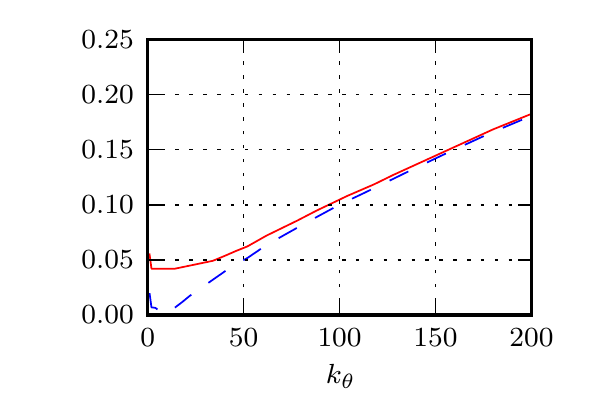}}
\caption{\label{thetaBiasMSEk} {\bf  $\hat{\theta}^{(k_\theta)}$ as a function of $k_\theta$: relative biases and MSEs.} Left: Log-Weibull distribution; Right: strict log-exponential-power law. From top to bottom: $\rho = 1.2, 2, 4$, $n = 1000$. Dashed lines: relative biases; Solid line: relative MSEs. These plot essentially shows the benefits of increasing $k_\theta$ from $1$ to $2$. They also indicate that $k_\theta$ can be further increased with benefits as long as $ k_\theta \ll n$. This is true for all $\rho$.}
\end{figure}

\begin{figure}
\centerline{}
\centerline{\includegraphics[width=60mm]{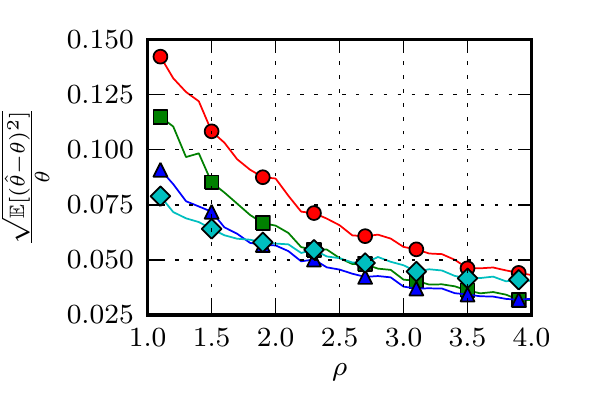} \hspace{1mm} \includegraphics[width=60mm]{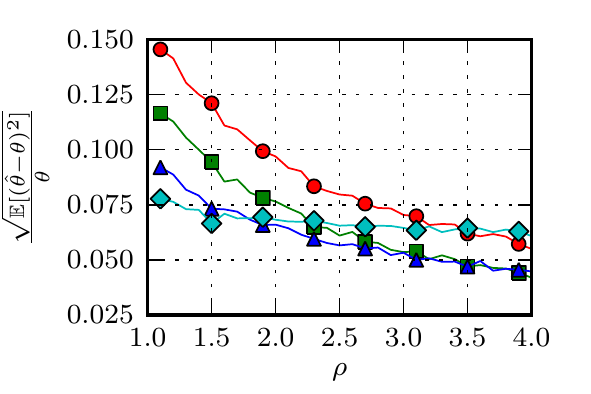} }

\centerline{\includegraphics[width=60mm]{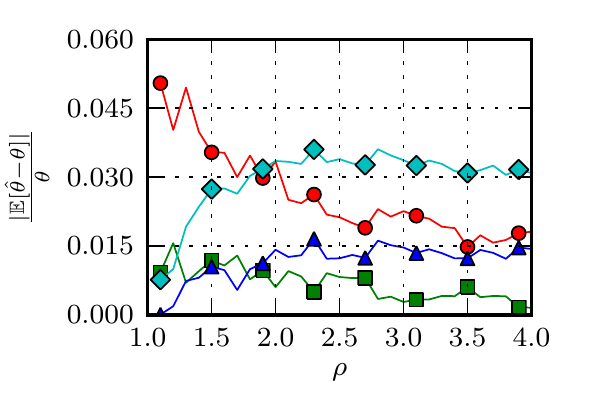} \hspace{1mm} \includegraphics[width=60mm]{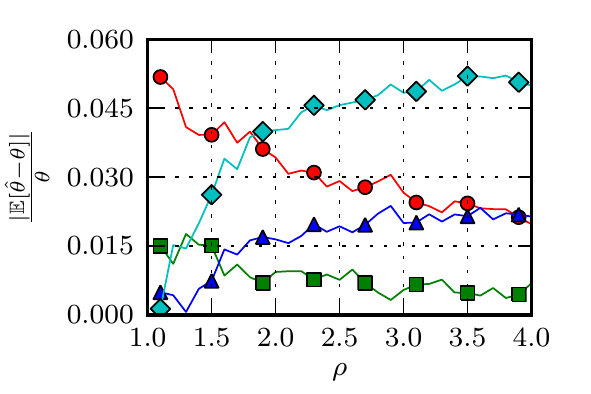} }

\caption{\label{thetaBiasMSErho} {\bf $\hat{\theta}^{(k_\theta)}$ as a function of $\rho$: relative biases and MSEs.} Left: Log-Weibull distribution; Right: strict log-exponential-power law. Top: relative MSEs; bottom: relative biases. $n = 1000$.  Red '$\circ$': $k_\theta= 1 $~; green '$\square$': $k_\theta=2$~; blue '$\triangle$': $k_\theta=25$~; cyan '$\Diamond$': $k_\theta=50$. These plots show again the benefit of increasing $k_\theta$ from $1 $ to $2$, and there exists an optimal $k_\theta$ that depends essentially on $n$ and weakly on $\rho$.}
\end{figure}

 \begin{figure}
\centerline{}
\centerline{\includegraphics[width=60mm]{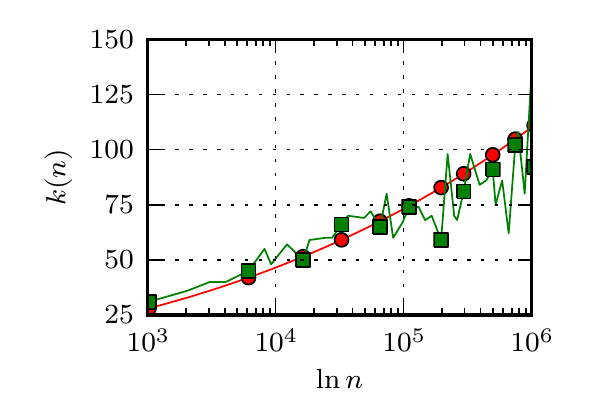} \hspace{1mm} \includegraphics[width=60mm]{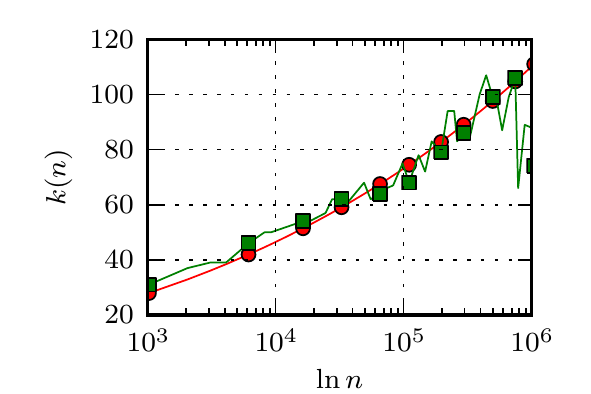} }

\centerline{\includegraphics[width=60mm]{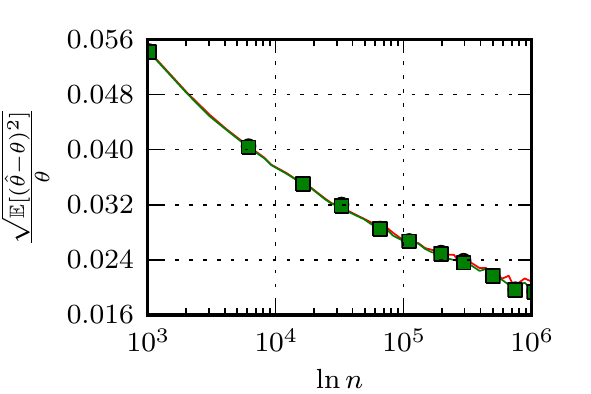} \hspace{1mm} \includegraphics[width=60mm]{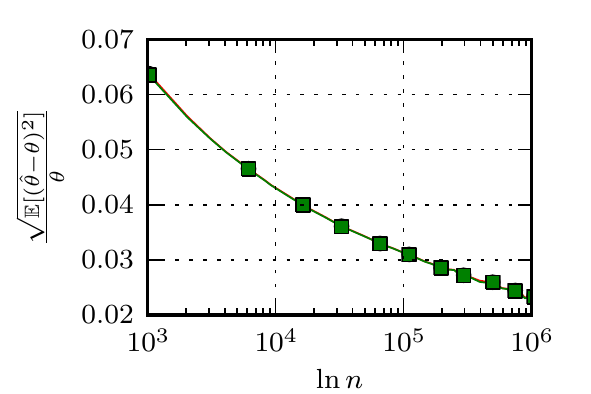} }

\caption{\label{thetaMSEn} {\bf  $\hat{\theta}^{(k_\theta)}$: $k_{\theta,opt}$ and MSEs as a function of $\ln n$. } Left: Log-Weibull distribution; Right: strict log-exponential-power law.
Top: $k_{\theta,opt}$, bottom: MSE. Red '$\circ$': $k_\theta=\exp( \sqrt{1.6 \ln n})$~; green '$\square$': $k_{\theta,opt}$ opt $= $ min of MSE. $\rho = 2$. These plots show the relevance of the proposed empirical formulae for $k_{\theta,opt}(n)$.}
\end{figure}

%
%
 
%
%
%

Numerical simulations are conducted using both categories of examples, log-Weibull and  strict log-exponential-power law distributions, described in Section \ref{sec:RV}, for different values of $\rho$, and different observation size $n$. 
For all simulations, $500$ independent realizations of $n$ \emph{i.i.d.} random variables are drawn. 

From these, the (relative) biases, variances and mean square errors (MSE) of $\theta^{(k)} $ are evaluated and reported in Fig. \ref{thetaBiasMSEk}.  
For both types of laws, it shows that, as expected from the theoretical analysis above, increasing $k_\theta$ from $1$ to $ 2$ yields a significant decrease in bias, and hence in MSE.
For $k_\theta = 2 $, $\alpha_1 \approx \alpha_2 \approx \frac{1}{2}$, hence showing that a simple arithmetic mean of the first and second largest observed values performs much better than the use of the sole maximum $M_{Y,n}$. 
The bias only slightly further decreases when $k_\theta$ is increased beyond $2$, before it finally increases again when $k_\theta$ becomes non negligible compared to $n$. 
Combining biases and variances, the MSEs show, for all $\rho$, the benefits of using $k_\theta \geq 2 $, yet small, $ k_\theta \ll n$, compared to $k_\theta=1$ or to $k_\theta$ too large (too close to $n$). 
Also, it can be noted that the closer $\rho $ to 1, the larger $ k_\theta $ needs to be to achieve the lowest MSEs. 

Furthermore, Fig. \reffig{thetaBiasMSErho} reports the relative MSEs (left) and biases (right) as a function of  $\rho$, for various choices $k_\theta$.
Besides confirming that the benefit from using $k_\theta =2$, or slightly larger, against $k_\theta =1$ is valid for all $\rho$, it also shows that biases and MSEs are roughly constant for $\rho \geq 2 $, while, in most cases, significantly increasing when $\rho \rightarrow 1^+$. 

To finish with, Fig. \reffig{thetaMSEn} compares $k_{\theta,opt}$ (i.e., $k_{\theta}$ that yields the lowest MSE) to tentative candidate functions for practical $k_\theta(n)$: $ k_\theta(n) = 10 \ln n - 40$, $ k_\theta(n) = \exp( \sqrt{1.6 \ln n})$. 
Though empirical, those formula hold, for both types of distributions, for a large range of values of $\rho$ ($\rho \in [1.5,3]$) and of $ n $, and can hence be used as rules of thumb for a practical use of the estimator. 

\subsection{Estimation of $\rho_l(y^\dagger(n))$}
\label{sec:rhoest}

\subsubsection{Definition}

For simplicity, $\rho_l(y^\dagger(n))$ is hereafter relabeled as $ \rho_E$. 
From its definition (cf.~Eq.~(\ref{eq:rhoe})),  $\rho_E$ can be interpreted as the local exponent governing the behavior of the function $h$ (cf.~Eq.~(\ref{eq:h})) around the largest available observation  $y^\dagger(n)$. 
This suggests the direct use of the empirical cumulative distribution function to estimate $\rho_E$. 
Indeed, for a sufficiently small interval around $y^\dagger$, one has $\partial_{\ln y} h(y^\dagger) = \rho_E$, which enables us to write: 
\begin{equation} 
  \ln (-\ln F_Y(y^\dagger)) \approx \rho_E \ln y^\dagger + \beta. 
\end{equation}
Considering the ordered list $Y_{i,n}$, 
the $(1-\frac{i}{n})$-th quantile can be approximated by $Y_{i,n}$, hence:
\begin{equation*}
 \ln \left(-\ln \frac{i}{n} \right) \approx \rho_E \ln Y_{i,n} + \beta.
 \end{equation*}
This suggests that $\rho_E$ can be estimated by a non weighted least square fit that involves $Y_{i,n}$ for $i \leq k_\rho$, where $ k_\rho$  is sufficiently small compared to $n$: 
\begin{equation} \hat{\rho}_E = \frac { C(\ln y, \ln(\ln n - \ln i) )} {C(\ln y,\ln y)}, 
\end{equation}
with $ \overline x= \frac{1}{k} \sum_{i=1}^{k_\rho} x_i $ and $ C(x,y)= \overline {xy} - \overline{x}  \,\overline{y}$.

\subsubsection{Performance: theoretical analysis}

Again, the choice of $k_\rho$ constitutes the critical issue that can be expressed in terms of the usual bias-variance trade-off: 
Because $\rho_E$ is a local exponent, that may converge only slowly to $\rho$, it needs to be estimated on narrow intervals around $y^\dagger $, hence the use of a small $k_\rho$; Conversely, large $k_\rho$ decrease the variance of the estimate at the price of a bias increase. 
Actually, a too large $k_\rho$ necessarily leads to the underestimation of $\rho_E$, as can be seen by evaluating the variations of $\rho_l(y)$ for each $Y_{1,n}$ and $Y_{k,n}$.
This motivates the choice of $k_\rho$ according to: 
\begin{equation}
\lim_{n  \rightarrow + \infty} \frac{k_\rho(n)}{n}=0.
\end{equation}

\subsubsection{Performance: numerical analysis}

 \begin{figure}
\centerline{}
\centerline{\includegraphics[width=60mm]{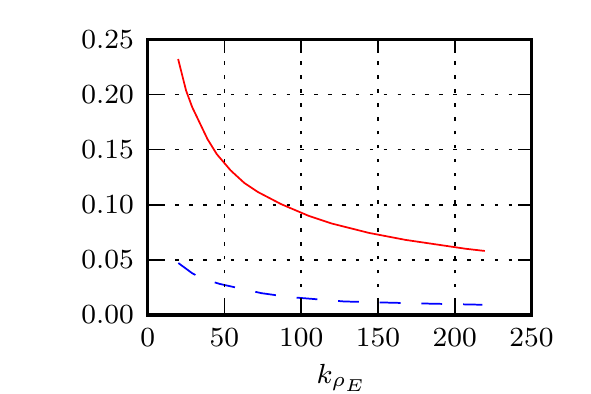} \hspace{1mm} \includegraphics[width=60mm]{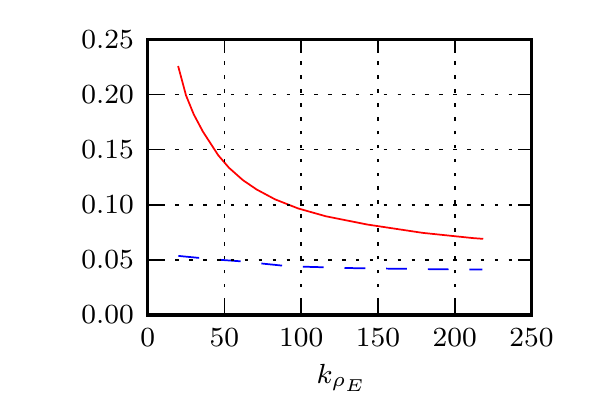}}
\centerline{\includegraphics[width=60mm]{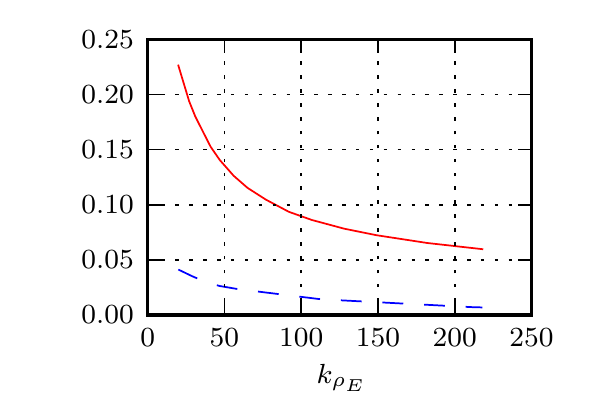} \hspace{1mm} \includegraphics[width=60mm]{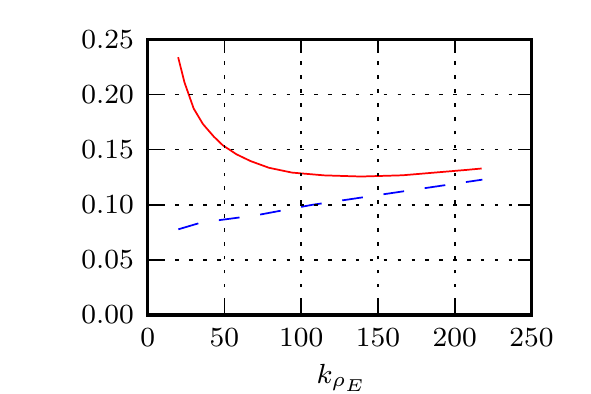}}
\centerline{\includegraphics[width=60mm]{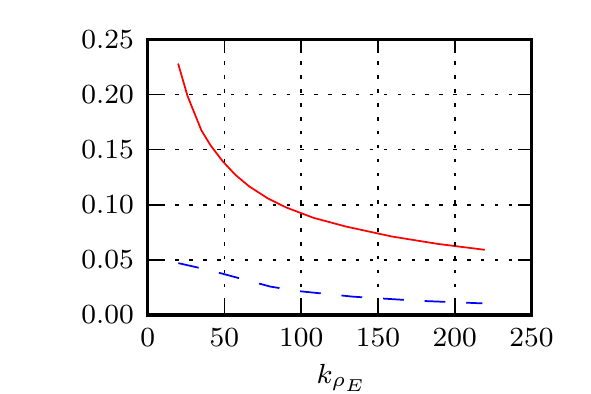} \hspace{1mm} \includegraphics[width=60mm]{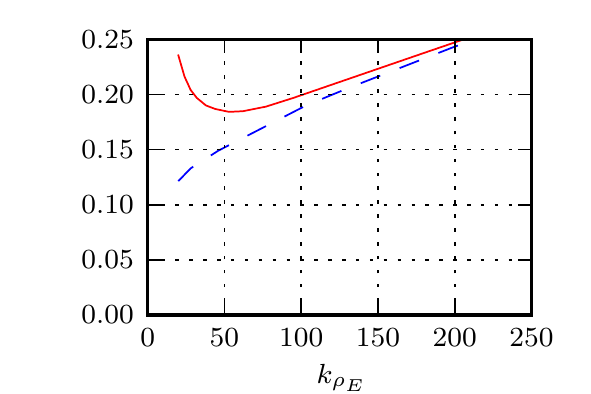}}

\caption{\label{rhoEBiasMSEk} {\bf $\hat{\rho}^{(k_\rho)}_E$ as a function of $k_\rho$: relative biases and MSEs.} Left: Log-Weibull distribution; Right: strict log-exponential-power. From top to bottom : $\rho = 1.2, 2, 4$, $n=1000$. Dashed blue lines: relative biases~; Solid red  lines: relative MSEs. For the Log-Weibull distribution, because $\rho_E(n) \equiv \rho$, increasing $k_\rho$ is beneficial as it only yields a reduction in variance. However, in general, there is an optimal $k_{\rho,opt}$ that depends on $n$ and $\rho$.}
\end{figure}

\begin{figure}
\centerline{\includegraphics[width=60mm]{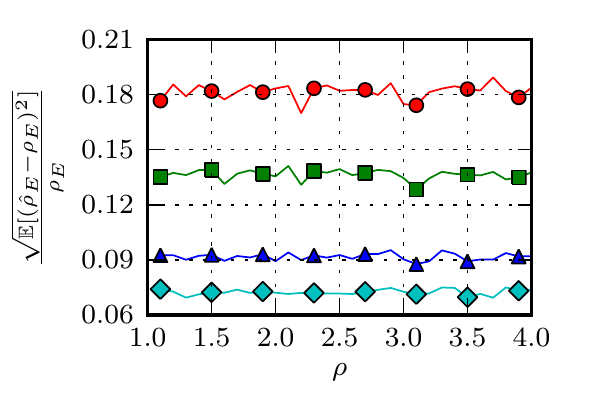} \hspace{1mm} \includegraphics[width=60mm]{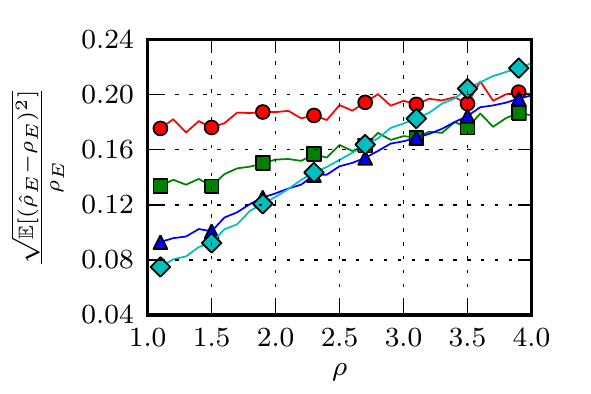} }
\centerline{\includegraphics[width=60mm]{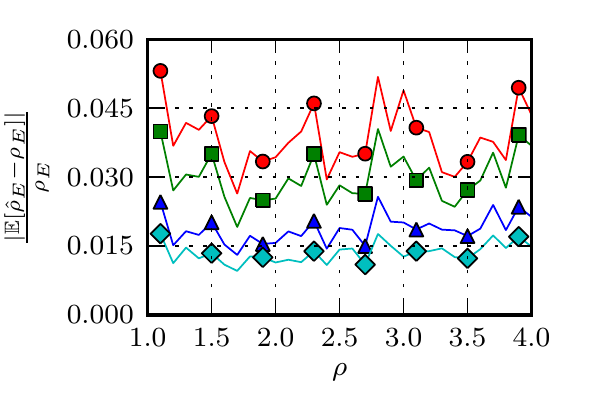} \hspace{1mm} \includegraphics[width=60mm]{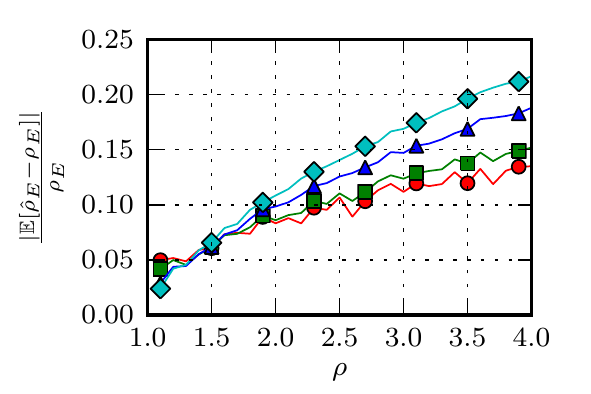} }
\caption{\label{rhoEBiasMSErho} {\bf  $\hat{\rho}^{(k_\rho)}_E$ as a function of $\rho$: relative biases and MSEs.} Left: Log-Weibull distribution; Right: strict log-exponential power. Red '$\circ$': $k_\rho=30 $~; green '$\square$': $k_\rho=50$~; blue '$\triangle$': $k_\rho=100$~; cyan '$\Diamond$': $k_\rho=200$.}
\end{figure}

 \begin{figure}
\centerline{}
\centerline{\includegraphics[width=60mm]{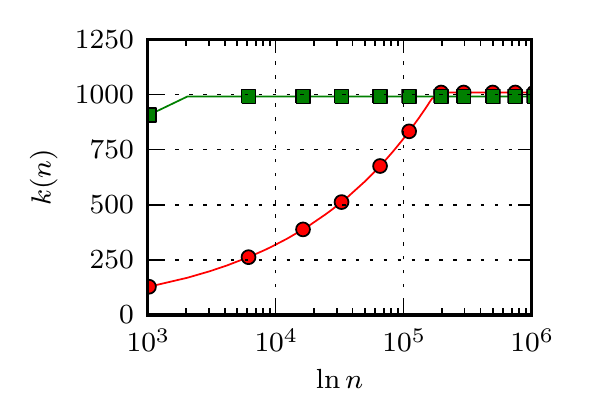} \hspace{1mm} \includegraphics[width=60mm]{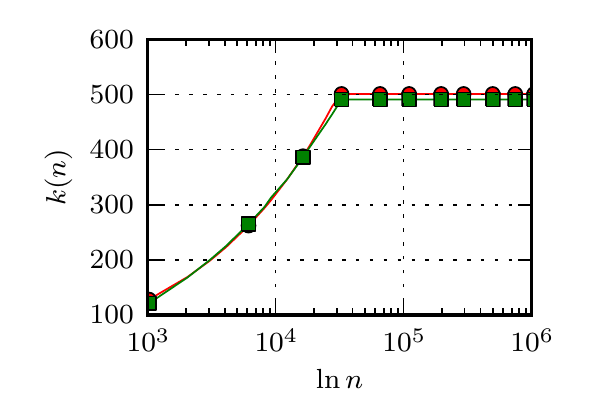} }
\centerline{\includegraphics[width=60mm]{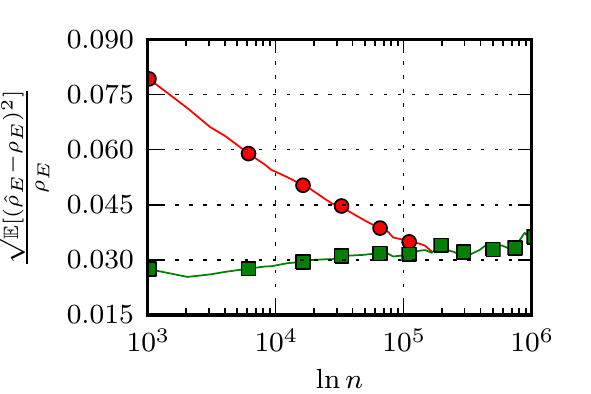} \hspace{1mm} \includegraphics[width=60mm]{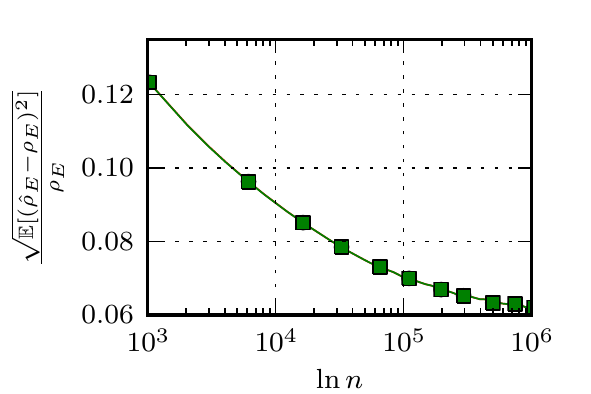} }

\caption{\label{rhoEMSEn}  {\bf $\hat{\rho}^{(k_\rho)}_E$: $k_{\rho,opt}$ and MSEs as a function of $\ln n$. } Left: Log-Weibull distribution; Right: strict log-exponential power.
Top: $k_{\rho,opt}$, bottom: MSE.  $\rho = 2$. Red '$\circ$' : $k_\rho(n)=8 n^{0.4}$~; green '$\square$': $k_{\rho,opt}$ opt $= $ min MSE for $k<1000$. $\rho = 2$. These plots show the relevance of the proposed empirical formulae for $k_{\rho,opt}(n)$.}
\end{figure}

The performance of $ \hat{\rho}_E$ as well as the trade-off framing the choice of $k_\rho(n)$ are studied numerically by means of Monte-Carlo simulations conducted as described in Section \ref{sec:thetanum}. 

Relative biases and MSEs are reported in Fig. \reffig{rhoEBiasMSEk}, as a function of $k_\rho$. 
It shows that for the log-Weibull distributions, bias is negligible and the larger $k_\rho$, the lower the MSEs. 
This stems from the fact that $\rho_E(n) \equiv \rho, \, \, \forall n $, 
which implies that $\hat{\rho}_E$ is unbiased whatever $k_\rho$, and hence that increasing $k_\rho$ decreases the variance at no cost. 
However, this is a particular case that departs from the general situation better accounted for by the second example consisting of strict log-exponential-power law distributions. 
On this second example, Fig. \reffig{rhoEBiasMSEk} indicates that the optimal choice of $k_\rho(n) $ depends on the actual value of $\rho$ and must be chosen larger for $\rho \rightarrow 1$. 

Relative biases and MSEs are reported in Fig.\reffig{rhoEBiasMSErho}, as a function of $\rho$. 
It shows that there is a uniform in $\rho$ benefit to increase $k_\rho$ for the log-Weibull distributions.
For the strict log-exponential-power law distributions, however, $k_{\rho,opt}$ is larger when $\rho \rightarrow 1$. 

To finish with, Fig.~\reffig{rhoEMSEn} compares $k_{\rho,opt}$ (i.e., $k$ that yields the lowest MSE) to tentative candidate functions for practical $k_\rho(n)$: $ k_\rho(n) = 50 \ln n  $, $ k_\rho(n) = 8 n^{1/3}$. 
Though empirical, those formula can be used as rules of thumb for a practical use of the estimator. 

\subsection{Estimation of $\qc$}
\label{sec:qcest}

\subsubsection{Definition}
From $ \hat \theta^{(k_\theta)}$ and $\hat \rho_E^{(k_\rho)}$ defined above, an estimator for  $\qc$ can be proposed as: 
\begin{equation} 
\hat{q}_c^{(k_\theta,k_\rho)}= \hat{\theta}^{(k_\theta)}  \hat{\rho}^{(k_\rho)}_E. 
\end{equation} 

\subsubsection{Performance: theoretical study}

To study the performance of $\hat{q}_c$, its bias  $\Bias{\hat {q}_c} =\E \hat q_c - \qc$ and variance $ \Var{\hat{q}_{c}}$ can be computed as
\begin{equation}
\Bias{\hat q_c}= 
\Cov{\hat{\theta}}{\hat{\rho}_E} 
+ \Bias{\hat \theta} (\rho_E + \Bias{\hat{\rho}_E})
+ \Bias{\hat{\rho}_E}(\theta+\Bias{\hat \theta} ) , 
\end{equation}
and
\begin{align*}
 \Var{\hat q_c}&= \Var{\hat{\theta}}\Var{\hat{\rho}_E}\\
&+ \Cov{\hat{\theta}^2}{\hat{\rho}_E^2}- \Cov{\hat{\theta}}{\hat{\rho}_E}^2\\
&+ \Var{\hat{\theta}}(\rho_E + \Bias{\rho_E})^2 + \Var{\hat{\rho}_E}(\theta+\Bias{\theta} )^2 - 2 \Cov{\hat{\theta}}{\hat{\rho}_E}  (\rho_E + \Bias{\rho_E}) (\theta+\Bias{\theta} ) \\ 
\end{align*}

Furthermore, assuming that the relative biases are negligible yields: 
\begin{equation}
 \frac{\Bias{\hat q_c}}{q_c } \approx   \frac{\Bias{\hat \theta}}{\theta}  + \frac{\Bias{\hat \rho_E}}{\rho_E}  + \frac{\Cov{\hat{\theta}}{\hat{\rho}_E}}{\rho_E \theta}.
 \end{equation} 
\begin{align*}
  \frac{\Var{\hat \qc}}{q^2_c}& \approx   \frac{\Var{\hat{\theta}}\Var{\hat{\rho}_E}}{(\rho_E \theta)^2} \\
&+ \frac{\Cov{\hat{\theta}^2}{\hat{\rho}_E^2}- \Cov{\hat{\theta}}{\hat{\rho}_E}^2}{(\rho_E \theta)^2}\\
&+ \frac{\Var{\hat{\theta}}}{\theta^2}  + \frac{\Var{\hat{\rho}_E}}{\rho_E^2} - 2 \frac{\Cov{\hat{\theta}}{\hat{\rho}_E}}{  \rho_E \theta}. \\ 
\end{align*}
These calculations show that the performance of $\hat{q}_c$ significantly depend on the covariance between the estimates $\hat{\theta}$ and $ \hat{\rho}$. This covariance is investigated by means of numerical simulations. 

\subsubsection{Performance: numerical analysis}

\begin{figure}  []
\centerline{
\includegraphics[width=60mm]{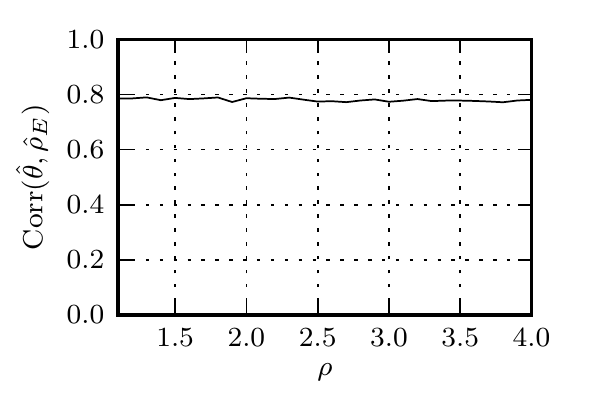} \includegraphics[width=60mm]{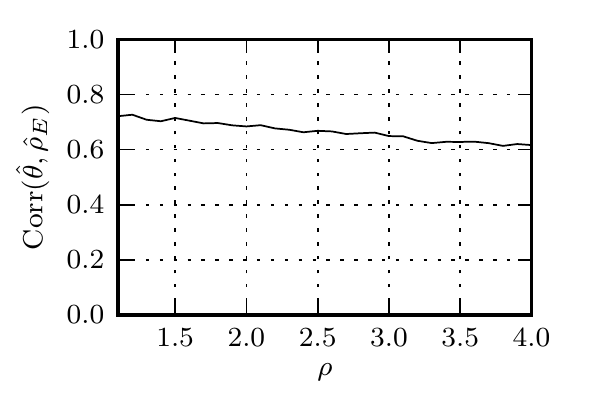}
}
\centerline{\includegraphics[width=60mm]{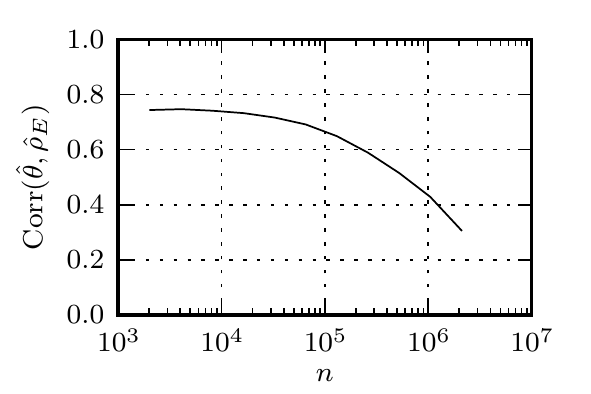} \includegraphics[width=60mm]{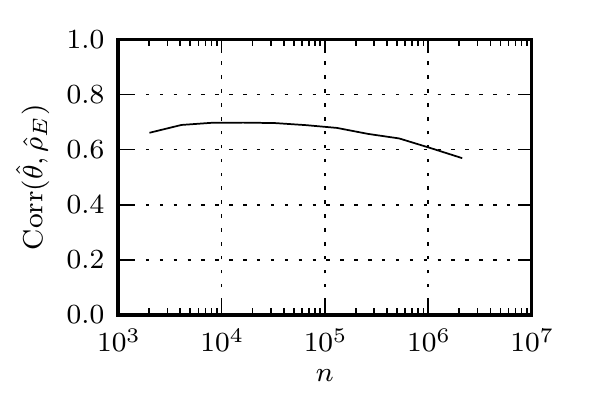}}
\caption{\label{fig:qc:correlationW} {\bf Correlation between $\hat \theta^{(k_{\theta,opt})}$ and $\hat \rho^{(k_{\rho,opt})}_E$. } Left: Log-Weibull distribution~; right: Strict log-exponential-power distribution. Top: as a function of $\rho$, $n=1000$. Bottom: as a function of $n$, $\rho = 2$.}
\end{figure}

\begin{figure}[]
\centerline{
\includegraphics[width=60mm]{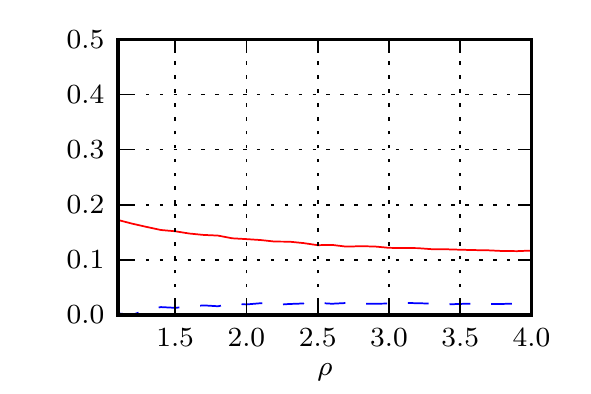} \includegraphics[width=60mm]{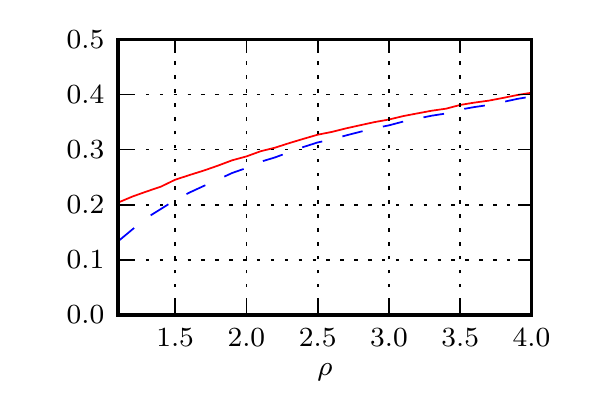}
}
\centerline{
\includegraphics[width=60mm]{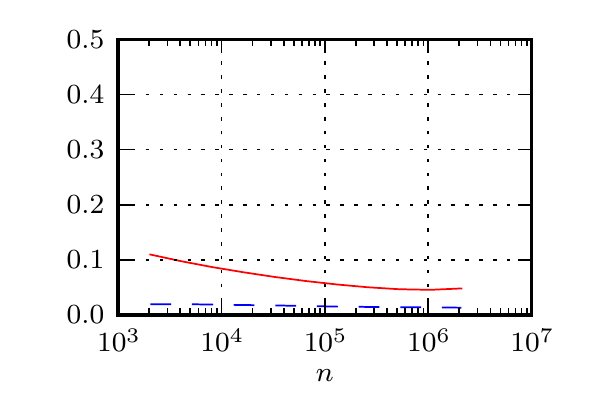} \includegraphics[width=60mm]{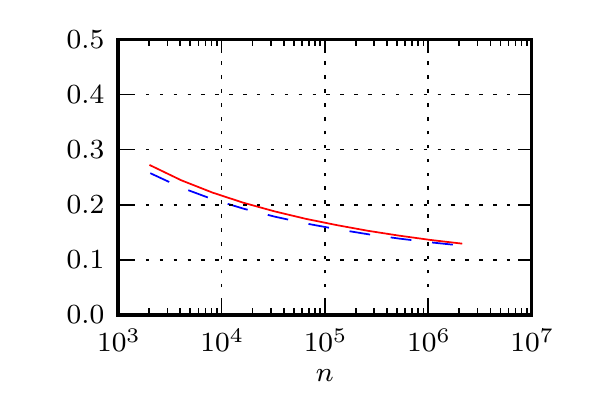}
}
\caption{\label{fig:qc:logwb} {\bf $\hat q_c$: relative biases (blue dashed line) and MSES(red solid line).}  Left: Log-Weibull distribution~; right: Strict log-exponential power distribution. Top: as a function of $\rho$, $n=1000$. Bottom: as a function of $n$, $\rho = 2$. $n=1000$}
\end{figure}

Monte-Carlo simulations, in the spirit of those presented in previous sections, with the empirically found best choices for $k_\theta(n) $ and $k_\rho(n)$, are conducted. 
Fig.~\reffig{fig:qc:correlationW} shows that $\Cov{\hat{\theta}}{\hat{\rho}_E}$ remains quite high and positive, whatever $\rho$ and $ n $, and hence significantly contributes to the bias and variance of $\hat \qc$. 
Furthermore, results reported in Fig.~\reffig{fig:qc:logwb} essentially show that the MSEs for $\hat q_c$ vary from $5 \% $ to $ 25 \%$ depending on the actual distribution of $Y$ and the value of $\rho$, with larger MSEs for $\rho \rightarrow 1$.  
Though these MSEs may seem high, $\hat q_c$ still provides practitioners with a practical procedure yielding a very satisfactory order of magnitude of the true $\qc$: 
In applications, one is rarely interested in the precise value of $q_c$, while the knowledge of a estimated order of magnitude often proves very important as soon as data interpretation is concerned. 



\section{Correlated time series}
\label{sec:cov}

So far, observation $\{X_1, \ldots, X_n \} $ consisting of independent samples only have been dealt with, when we are mostly interested in stationary time series analysis, where data are likely to be dependent. 
The goals of the present section is therefore to analyze the impact of dependence amongst observations on the definition and estimation of $q_c$ as well as to propose a simple and elementary modification of $\hat{q}_c$ as described above to accommodate dependencies. 

\subsection{Theoretical analysis}

In Section \ref{sec:analyse} above, the critical moment arose from a balance between two competing factors: dominant moment contribution versus finite size effect.
While the definition and evaluation of the former depend only on the marginal distribution of the data and are not affected by dependencies, the latter is far more impacted by correlations. 
Indeed, a finite size accessible range has been defined by  
 \begin{equation*} 
\Prob{\forall k=1, \ldots, n, Y_k < y^\dagger} = \frac{1}{e}, 
\end{equation*}    
For independent data, the previous equation simplifies to $h(y^\dagger)=\ln n$.
For correlated time series, it is natural to expect that $y^\dagger_T$ can be defined via an effective number of sample $n^*_T$, as:
 \begin{equation*} 
h(y^\dagger_T)=\ln n^*_T. 
\end{equation*} 
From this modification, the analysis conducted in Section \ref{sec:analyse} leads to define a critical moment under dependencies $q_c^T$ as:  
 \begin{equation} 
 \label{eq:qc:corr}
q_c^T(n)= q_c(n^*_T). 
\end{equation} 
As for the independent case, $q_c^T$ can be split into two terms,
 \begin{equation} 
 \label{eq:qc:corr;thetaRho}
q_c^T(n)= \rho_E^T(n) \theta^T(n). 
\end{equation} 
Reinjecting $y^\dagger_T$, into the definition of $\rho_E^T$ and $\theta^T$ leads to 
  \begin{equation} 
\rho_E^T(n)= \rho_E(n^*_T), \quad
\theta^T(n)=\theta(n^*_T). 
\end{equation} 

Let $\tau$ denote the correlation length of the analyzed time series, defined as $ \tau= \int_{0}^{+\infty} t C(t) dt / \int_{0}^{+\infty} C(t) dt .$
When $\tau$ is well defined,  it is natural to assume that $n^*$ essentially behaves as $n^* \propto \frac{n}{\tau}$, with the obvious requirement that
$ \lim_{\tau\rightarrow 0} n^*(n,\tau)=n$, suggesting the phenomenological form: 
\begin{equation}
\label{eq:nstar}
n^* = \frac{n}{1+\kappa \tau},
\end{equation}
with $\kappa$ a positive constant to be determined. 
This leads to: 
\begin{equation}
\label{eq:qctau}
q_c^T(n)=q_c(\tau, n)= q_c(n/(1+\kappa \tau)).
\end{equation}

\subsection{Impact of dependencies on $\hat q_c$} 

To study the impact of dependence on $\hat q_c$, the estimator is applied to a large number of independent realizations of stationary time series, of length $n$, and with a priori and jointly prescribed (log-Weibul or strict log-exponential power law) marginal distributions and covariance function $C(t)$. These are synthesized numerically using the Hermite expansion and Circulant Matrice Embedding based techniques described in \cite{hpa11a,hpa11b}. To understand the impact of dependencies, let us analyze the case where $C(t)$ takes the simple form:  
\begin{equation} 
\mathrm{C}(t)= e^{-\frac{|t|}{\tau}}, 
\end{equation}
with $\tau $ a tunable constant, obviously corresponding to the correlation length. 

\begin{figure}
\centerline{ \includegraphics[width=70mm]{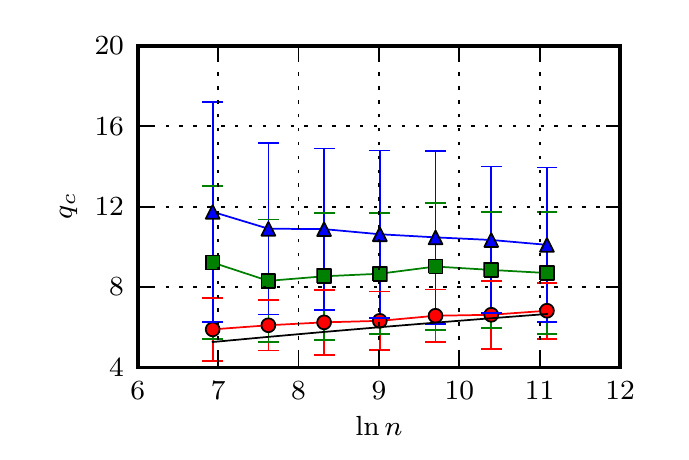} }
\caption{\label{fig:qc:uncorrected}{\bf $\hat q^{(k_\theta,k_\rho)}_c$: Impact of dependencies.} Applied to correlated observations, $\hat q_c$ is impaired by a significant structural bias that cannot be accounted for by a rescaling (of the form suggested by Eq.  \protect (\ref{eq:qctau})). $(k_\theta,k_\rho) = (1,100)$. Black line: true $q_c$, $\hat q^{(k_\theta,k_\rho)}_c$, for  $\tau=10$ (red circles), $\tau=50$ (green squares) $\tau=100$ (blue triangles).} 
\end{figure}

Results, reported in Fig.~\reffig{fig:qc:uncorrected}, show first a strong discrepancy between $\hat q_c$ and $q_c(\tau,n)$ and second and mostly that no rescaling operation of the form suggested by Eq. (\ref{eq:qctau}) can account for it. 
Indeed, a rescaling as in Eq. (\ref{eq:nstar}) amounts only to a translation of  $\hat{q}_c$ along the horizontal axis. 
Moreover, assuming $n^*<n$  restricts to left translations only, which makes it impossible to superimpose all $\hat{q}_c(\tau,\n)$ curves by rescaling: 
 $\hat{q}_c$ does not converge towards $q_c(\tau,n)$.  

To further analyse the cause for this discrepancy, let us analyse separately $\hat \theta^{(k_\theta)}$ and $\hat \rho^{(k_\rho)}_E$. 

\begin{figure} [!h]
\centerline{\includegraphics[width=70mm]{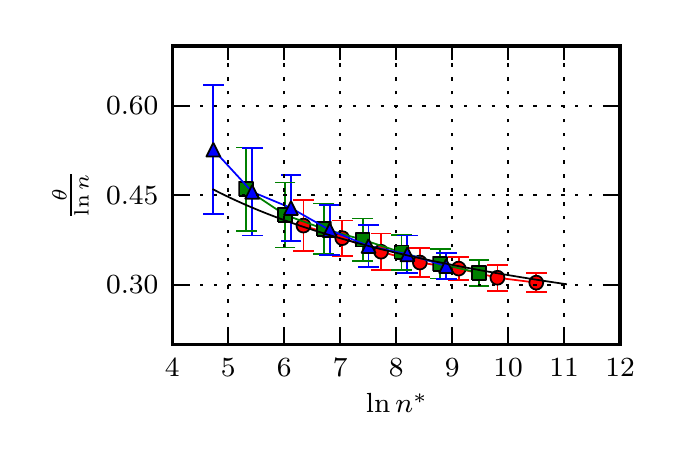}}
\caption{\label{fig:theta:scaling}{\bf $\hat{\theta}$: Impact of dependencies.} Rescaling of $\frac{\theta}{\ln n}$ as a function of $n^* = n/(1 + \kappa \tau)$,  for $\tau \in \{10,50,100\}$ (red circles $\tau=10$, green squares $\tau=50$ , blue triangles $\tau=100$). The black line correspond to the theoretical value $\theta(n)/\ln n$. }
\end{figure}

Fig.~\reffig{fig:theta:scaling} obtained from the same Monte-Carlo simulations, with $k_\theta = 1$, indicates that  a $\kappa (\approx 0.08)$ can be found such that $\Esp {\frac{\hat{\theta}} {\ln n}  } \approx \frac{\theta(n^*)}{\ln n^*} $ (quite logically $\hat{\theta} = \ln n^* / \Omega $ is to be used instead of $\hat{\theta} = \ln n / \Omega $, as $\ln n$ in the definition of $\theta$ does not scale with $\tau$.)
This means that $\Omega_1$ exhibits the expected behaviour, $\Esp{ \Omega_1(n,\tau)}=\Esp{\Omega_1(n^*)}$. 
However,  increasing $k_\theta$ from $1$ to $2$ and above does not, as in the independent case, brings a reduction in bias, but instead significantly increases it. 
This can easily be understood as very likely the second largest value is located extremely close to the first largest one and highly correlated to it, hence not bringing any bias reduction. 
Instead, this correlation prevents the convergence of the order statistics to its asymptotic limit and hence impairs the consistency of $\Omega$. 

\begin{figure} [!h]
\centerline{\includegraphics[width=70mm]{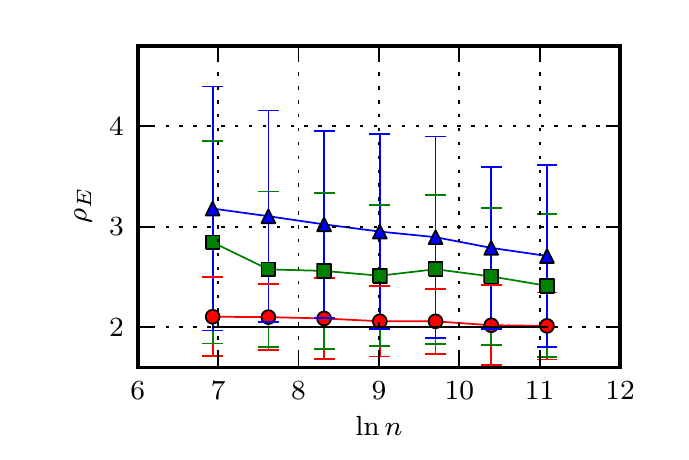}}
\caption{\label{fig:rho:alpha:0} {\bf $\hat{\rho}_E$: Impact of dependencies.} Black line: theoretical value, red circles $\tau=10$, green squares $\tau=50$ , blue triangles $\tau=100$. }
\end{figure}

This significant correlation of the order statistics also explains the poor performance of $\hat{\rho}_E$.
 As illustrated in Fig.~\reffig{fig:rho:alpha:0},  $\hat{\rho}_E$ does not converge towards $\rho_E(\tau,n)$ and no rescaling can account for this discrepancy. 
 This can be interpreted as follows:  $\hat{\rho}_E$ strongly relies on the convergence of the $k$th order statistics towards the $(1-k/n)$-quantile. Correlations of the sample generates correlations amongst this order statistics because correlated samples tends to regroups together in clusters in terms of ranks. So the correlation of the samples is transferred to the order statistics, which dramatically hinders the convergence of the order statistics towards the quantile.

Therefore, applied to correlated stationary times series,  $\hat{q}_c$, as defined in Section \ref{sec:estim}, cannot converge towards $q_c^T$. 
However, a simple modification described below is shown to significantly account for correlations: It consists of imposing that the largest value has been found, the second largest value used must be picked \emph{far enough} from the first largest one, the third largest must be picked far  \emph{far enough} from the two first ones, and so on. 
The modified $\hat{q}_c^T$ hence essentially consists of the original $\hat{q}_c$, applied to these \emph{far enough} order statistics. 


\subsection{Modified estimators} 

An efficient way of implementing the idea described above is to construct the set of effective independent samples step-by-step. 
Consider the original sample $\mathcal{P}_1=X_1,\dots, X_n$. 
First, the maximum, $X_{k_1}=\max \{\mathcal{P}_1\}$ is selected. 
Then, all samples located \emph{close} to $X_{k_1}$ in time are eliminated:  
using the distance 
\[ d_\beta(X_i,X_j) = \max(|j-i|,\beta \mathrm{Card} \{k,\min(X_i,X_j)<X_k<\max(X_j,X_i)\}). \]
the set $\mathcal{P}_1$ is sieved by eliminating samples belonging to the sphere centered at $X_{k_1} $ and of radius $s$. 
A new set $\mathcal{P}_2=\mathcal{P}_1 \setminus B_{d_\beta}(X_{k_1},s)$ is defined and the procedure is repeated until  a set $\mathcal{P}_e=\{X_{k_1},...,X_{k_l}\}$ of effective independent points is obtained. 
The modified $\hat{\theta}^T$ and $\hat{\rho}_E^T$ are defined by applying the original $\hat{\theta}$ or $\hat{\rho}_E$ to $\mathcal{P}_e$. 

\begin{figure}
 \centerline{\includegraphics[width=70mm]{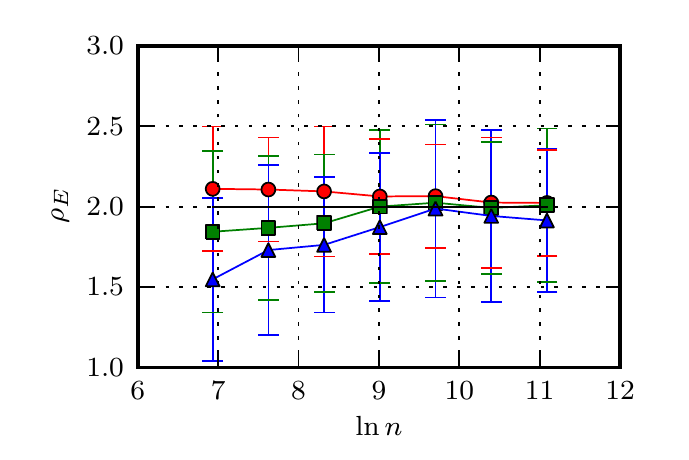} } 
\caption{\label{fig:corr:rho} {\bf Modified estimator $\hat{\rho}_E^T$.} $k_\rho=100$. Black line: theoretical value, red circles $\tau=10$, green squares $\tau=50$ , blue triangles $\tau=100$. } 
\end{figure}

Monte-Carlo simulations conducted as above suggest that the threshold should take the following form: 
\begin{equation} 
s = \alpha \tau,
 \end{equation}   
with the choice $\alpha \approx 0.01$, as illustrated in Fig.~\reffig{fig:corr:rho}. 

\begin{figure}
\centerline{ \includegraphics[width=70mm]{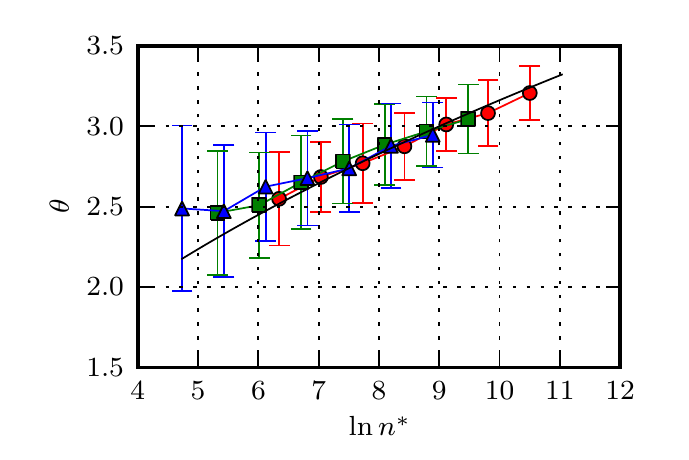} }
\caption{ \label{fig:theta_corr_rescaled} {\bf Estimation of $\theta(n^*)$ with the corrected estimator $\hat{\theta}^T$.} $k_\theta=10$. Black line : theoretical value, red circles $\tau=10$, green squares $\tau=50$, blue triangles $\tau=100$.} 
\end{figure}

\begin{figure}
\centerline{\includegraphics[width=70mm]{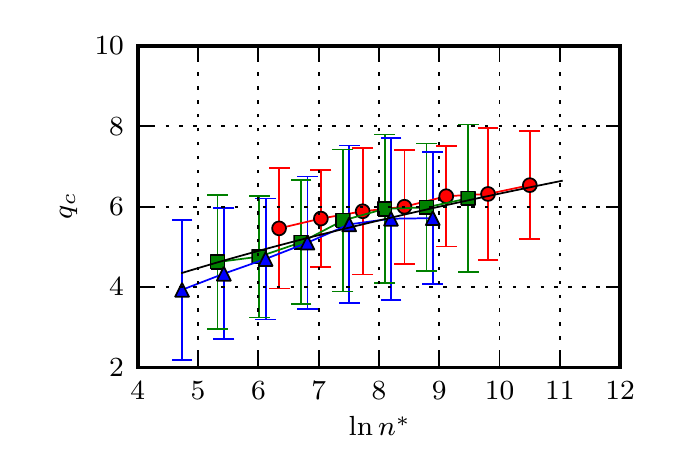} }
\caption{  \label{fig:qc_corr_rescaled}  {\bf Estimation of $q_c(n^*)$ with corrected estimator.} Black line: theoretical value, red circles $\tau=10$, green squares $\tau=50$, blue triangles $\tau=100$. } 
\end{figure}

These same Monte-Carlo simulations also show that the estimation of $\theta$ once corrected (by assuming the value $n^*$ is known) leads to the expected $\theta(n^*)$ (cf. Fig.~\reffig{fig:theta_corr_rescaled}) and that $q_c$ can thus be very satisfactorily estimated (cf. Fig.~\reffig{fig:qc_corr_rescaled}).

\begin{figure}
\centerline{\includegraphics[width=70mm]{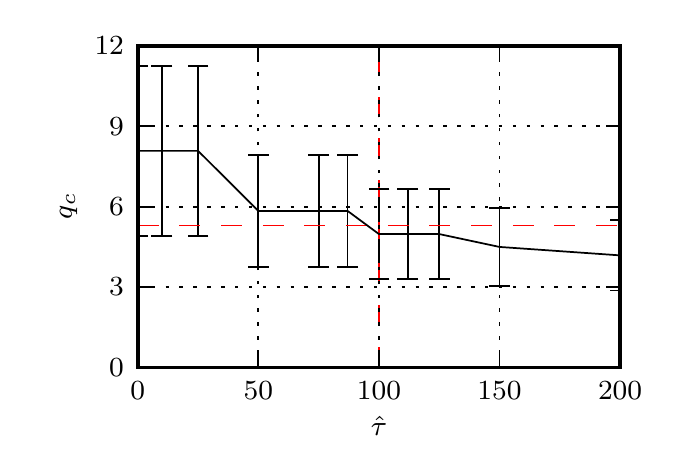} }
\caption{\label{fig:tau_influence} {\bf $\hat q^T_c $: Influence of the choice of $\tau$.} Solid black lines show the average value and standard deviations of  $\hat q^T_c $ using various values of $\tau$. The vertical  dashed line corresponds to the actual correlation length of the data, $\tau=100$. 
The horizontal dashed line indicates the theoretical value $q_c(n,\tau)$.} 
\end{figure}

This very good performance however involved the a priori knowledge of $\tau$, and Monte-Carlo simulations provided accurate estimates of the parameters $\alpha$ and $\kappa$.
These parameters would obviously not be known on actual data and would need to be estimated. 
Such estimations are beyond the scope of the present contribution.
They could be based on spectrum estimation or other classical techniques.
This is under current investigation and will be the subject of another contribution.
However, to conclude this study, Fig.~\reffig{fig:tau_influence} reports the expectations and standard deviations of $\hat q^T_c$, applied to log-normal exponentially correlated data with $\tau = 100$, when the value of $\tau $ actually used on $\hat q^T_c$ is, in purpose, incorrectly specified. 
Satisfactorily, Fig.~\reffig{fig:tau_influence} shows a weak sensitivity of $\hat q^T_c$ with the value of $\tau$ actually used, notably when it is larger than the actual correlation length. 
This means that, in practice, a rough estimation of the true correlation length might be sufficient to obtain satisfactory performance of $\hat q^T_c$. 


\section{Conclusions}

In the present contribution, the question of defining and estimating a critical order $q_c$ up to which the moments can be correctly estimated  from a sample observation of given size has been addressed. 
More specifically, the case where the observations stem from distributions that have finite moments of all orders yet significantly heavier tails than Gaussian random variables (such as e.g., the log-normal distribution) has been considered as it arises in many modern applications, such as e.g., the Random Energy Model in statistical physics or Multiplicative Cascades in Hydrodynamic Turbulence. 

It has been shown that this critical order can be given a formal quantitative definition, combining dominant contributions and finite size effects. 
This definition clearly highlights the impact of the parameter $\rho$ quantifying the heavytailness of the distribution on the growth of the critical order $q_c$ as a function of the size of the sample (the more standard case of variables having tail of the order of that of Gaussian random variables or a power-law tail being retrieved by taking the limit $\rho \rightarrow + \infty$ and $\rho \rightarrow 1$). 
Monte-Carlo simulations validates the practical relevance of this quantitative prediction. 
An estimator that can actually be applied to a single realization of data has been defined, and analyzed both theoretically and by means of Monte Carlo simulations, in the context of i.i.d. observations. 
This procedure provides practitioners with a tool that can practically be applied to a single observation of a given actual data.
{\sc Matlab} procedures implementing this estimation procedure are available upon request.

The question of dependence amongst observations has then been addressed and shown to significantly decrease estimation performance. 
However, a simple modification of the estimation procedure has been devised to account for dependence.
Its effectiveness is validated by Monte-Carlo simulations. 
The use of this procedure on actual data only requires a rough estimation of a typical correlation length amongst observations, which can be achieved using classical signal processing tools (such as spectrum estimation). 

The definition and assessment of an original procedure systematically integrating the estimation of the typical correlation length within that of $q_c$ is under current investigation and will be the subject of a future contribution.

An analysis, elaborated along the same line, has recently been proposed to define the critical moment for another class of variables: that corresponding to the increments of multifractal processes, that differ (significantly) from the class considered here in having finite moments only over a finite range of orders as well as long range correlations.
The reader is referred to \cite{Angeletti2010} for comparisons enabling insightful understanding of the mechanisms at work in defining a critical order.

\bibliography{biblio,math_analysis,math_proba,self} 

\begin{appendix}

\section{Saddle-point evaluation of $\E X^q$}
\label{app:SaddlePoint}

The saddle-point method consists in a classical technique enabling the evaluation of  the asymptotic behaviour of integrals of the form 
$ \int_{-\infty}^{+\infty} g(x) e^{t f(x)} dx $ for large $t$ \cite{Royden,DeBruijnAsymptotic}.  
More precisely, when $x^*$ is the only maximum of $f$, the saddle-point evaluation reads : 
\begin{equation} 
\label{eq:classicalSaddle}
\int_{-\infty}^{+\infty} g(x) e^{t f(x)} dx \simInf{t} g(x^*) e^{t f(x^*)}, \, \, \, \makebox{ with }  f'(x^*) = 0. 
\end{equation}

This technique is used here to evaluate the asymptotic behaviour of $\Esp {X^q}$ when $q\rightarrow +\infty$. 
Eq. (\ref{eq:EXqbis}) leads to define the logarithmic density $ \psi(y)=h(y)+\ln h'(y) $. 
The maximum of the integrand in Eq. (\ref{eq:EXqbis}) is hence located at $y^*$, such that $\psi'(y^*)=q$. 
The change of variable $y= \nu y^*$ enables to rewrite Eq. (\ref{eq:EXqbis}) as: 
\begin{equation}
\Esp{X^q}=y^* \int_0^{+\infty} e^{q y^* \nu - \psi(\nu y^*) } d\nu.
\end{equation}
Moreover, Conditions \ref{eq:hconda} and \ref{eq:hcondb} imply that  $\psi(y)\simInf{y} h(y) $ and $ \psi'(y)\sim_{y\rightarrow +\infty} h'(y) $. 
Combined together, these results yield: 
\begin{equation*}
q y^* t - \psi(t y^*)= \psi(y^*) ( \rho \nu - \nu^\rho + \varepsilon(y^*,\nu)) , 
\end{equation*}
with $\forall \nu, \lim_{y^*\rightarrow +\infty} \varepsilon(y^*,\nu)=0$; and hence
\begin{equation}
\Esp{X^q}=y^* \int_0^{+\infty} e^{ \psi(y^*) ( \rho \nu - \nu^\rho +\varepsilon(y^*,\nu))  } d\nu.
\end{equation}
The above equation can also be rewritten as
\begin{equation*}
\Esp{X^q}=y^* e^{ q y* -\psi(y^*) } \int_0^{+\infty} e^{ \psi(y^*) ( 1- \rho + \rho \nu - \nu^\rho +\varepsilon_2(y^*,\nu)) } d\nu.
\end{equation*}
with $\varepsilon_2(y^*,\nu)=\varepsilon(y^*,\nu)-\varepsilon(y^*,1)$. 
For $\nu\neq 1$, $1-\rho + \rho \nu - \nu^\rho +\varepsilon_2(y^*,\nu))<0$ thus only the neighbourhood of $\nu=1$  contributes significantly to the integral when $y\rightarrow+\infty$. Moreover, for $\nu$ close to $1$
\begin{equation*}
1-\frac{y^*\psi'(y^*)}{\psi(y^*)} + \rho \nu - \nu^\rho +\varepsilon_2(y^*,\nu)) \approx -y^{*2} \frac{\psi''(y^*)}{\psi(y^*)} \frac{\nu^2}{2} +o(\nu^2)
\end{equation*}
Injecting this approximation into the integral above leads to a Gaussian integral
\begin{equation*}
\int_0^{+\infty} e^{  -y^{*2} \psi''(y^*) \nu^2/2  } d\nu= 
\sqrt{\frac{2 \pi}{y^{*2} \psi''(y^*)}}, 
\end{equation*}
and  hence to Eq.(\ref{eq:momsaddlepoint}) :
\begin{equation}
\Esp{X^q} \sim e^{ q y^* -\psi(y^*) } \sqrt{\frac{2 \pi}{\psi''(y^*)} }   .
\end{equation}


%
%

\section{Truncated moments and sample moment estimator}
\label{app:LogConv}
In this section, the almost sure convergence of $f(n,q)=\ln S(n,q) / \ln n$ towards 
$\lim_n \ln M_T(n,q) / \ln n $ is proven.

The most fundamental point is to make precise the role play by $y^\dagger_\tau$ as a frontier. 
By definition,
\begin{equation*}
 \Prob{ Y_1 , \dots , Y_n > y^\dagger} = 1-e^{-\tau} \underset{\tau\rightarrow 0}\sim \tau 
 \end{equation*}
The Borel-Cantelli lemma \cite{Billingsley} states that if a sequence of random events $A_n$ satisfies $ \sum P(A_n) < +\infty$, then the events $A_k$ only occurs a finite number of times. 

Here, a natural choice of events would be $A_n=(\exists i<n,\, Y_i > y^{\dag}_\tau{n})$.  If $ \sum_n \tau(n) < +\infty $ then for sufficiently large $n$, all the $Y_i$ are almost surely smaller than $y^\dag_k$.
However, to preserve the property $ y^\dagger_\tau \rightarrow y^\dagger_1$, it is also required that $ \ln \tau \ll \ln n$. Unfortunately, these two conditions are incompatible.
 
To answer this dilemma, let us consider a subsequence $n(r)$ of samples (for instance $n(r)=2^r$ ) and events $A_r=(\exists i<n(r),\, Y_i > y^{\dag}_\tau{n(r)})$. 
Then, the condition in the Borel-Cantelli lemma becomes :
 \begin{equation}
 \sum_r \tau(n(r)) < +\infty 
 \end{equation}
One should note that this condition, with $n(r)=2^r$, is compatible with $\logLt \tau$ if one choose $ \tau= 1 / \ln n$. 
Moreover, if the value of $\ln S/ \ln n$ at points $n(r)$ is known, the behaviour of $f$ can be inferred for all points.
Considering two points $k<l$, one has:
\begin{equation*}
f(k,q(k)) < f(l,q(l)) \frac{\ln l}{\ln k} + 1 - \frac{\ln l}{\ln k} 
\end{equation*} 
With $a=n(r)<k<b=n(r+1)$, the inequality below holds: 
\begin{equation}
f(a,q(a)) \frac{\ln a}{\ln b} + 1 - \frac{\ln a}{\ln b} < f(k,q(k)) < f(b,q(b)) \frac{\ln b}{\ln a} + 1 - \frac{\ln b}{\ln a}.
\end{equation}
With $\ln(n(r+1)) / \ln(r) \rightarrow 1$, the above inequality describes precisely $f$, in the limit $n \rightarrow +\infty$. 
This equivalently amounts to defining $n(r) $ as 
\begin{equation*}
n(r) = \exp( e ^ {\nu(r)} ),\, \text{with} \lim_{r \rightarrow +\infty} \nu(r+1)-\nu(r)=0
\end{equation*}
A simple solution would be to choose $\nu(r)=\sqrt r$ and $\tau(r)=\frac{1}{r^2}$, then $ n(r)= \exp(e^{\sqrt r})$ and 
 \begin{equation*}
 \tau(n)= \left( \frac{1}{\ln \ln n} \right)^4.
\end{equation*}
Therefore, if points $n(r)$ only are observed, for $r$ sufficiently large, almost surely, all the $Y_i$ are smaller than $y^\dagger_\tau$. 

Denoting $\chi_I$ the characteristic function of the set $I$,
\begin{equation}
\Char{I}{x} = \begin{cases} 1 & \text{if } x\in I \\ 0 & \text{otherwise} \end{cases},
\end{equation}
one can write
\begin{equation*}
S(n,q) \as{=} \frac{1}{n} \sum Y_i \Char{(-\infty,y^\dagger]}{Y_i} 
\end{equation*}
This sum can be further evaluated by spliting the remaining interval $(-\infty,y^\dagger_\tau]$ in $2 R(n) +1$ sub-intervals with 
\[ a_k=k  \frac{y^\dagger}{R(n)} \]
\[ I_k= ( a_k, a_{k+1}  ], k \in -R,\dots, R-1. \] 
\[ I_{-R(n)-1}=I_{-\infty}= ( -\infty, a_{-R}  ]  \] 
For reasons made explicit later, we choose
\begin{equation}
R(n)= \frac{\ln n}{\ln \ln n}.
\end{equation}
Because $y^\dagger =\mathcal{L}(\ln n) (\ln n)^{1/\rho}$, the length of the sub-intervals is converging towards $0$.
$S(n,q)$ can hence be rewritten as
\[ S(n,q) \as{=} \sum_{i=1}^{n} \sum_{k=-R(n)-1}^{R(n)-1} e^{q Y_i} \Char{I_k}{Y_i}. \]
It is possible count the number of point in each interval: 
\[ N_k = \sum_{i=1}^{n} \Char{I_k}{Y_i} . \]
Under certain conditions $N_k$ should not differ too much from
 \begin{equation}
 \E N_k = n \E \chi_{I_k}
\end{equation}
More precisely, Markov inequality \cite{Billingsley} provides us with an almost sure upper bound:
\[ \Prob{\forall k, N_k > \alpha E N_k} < \frac{2 R(n)+1}{\alpha(n)} \]
To use Borell-Cantelli lemma, it is required that $\sum_r \frac{2 R(n(r))+1}{\alpha(n(r))} < +\infty$. 
The slower $\alpha$ grows, the more information can be obtained from Borel-Cantelli lemma. 
A valid choice here implies $\logLt \alpha$. 
Due to the sparsity of $n(r)$, 
 \begin{equation}
\alpha(n)= R(n) \ln n,
\end{equation}
 is a compatible choice.
 This yields:
\begin{equation}
 S(n,q) \as{<}  \sum_{k=-R_n-1}^{R(n)-1}  \alpha(n) \Esp {\CharF{I_k}} e^{q a_{k+1}}
\end{equation}
Furthermore, let $\psi$ denote the logarithmic density (as defined in \ref{app:SaddlePoint}),
\[ 
\Esp {\CharF{I_k}} < 2 (y^\dagger/R) \exp(-\psi(a_k)).
 \]
The major contribution to the sum comes from the interval $I_m$ containing $y^*$ (defined as  $\psi'(y^*)=q$). 
This leads to a rough upper bound,
\begin{equation}
S(n,q) \as{<} \alpha(n) \left( y^\dagger_\tau e^{q a_{m+1}-\psi(a_m)} + \Esp{\chi_{I_{-\infty}}} e^{- q y^\dagger} \right) .  
\end{equation}
Most of the terms in this equation are logarithmically negligible compared to $\ln n$. Therefore,
\begin{equation}
 \label{eq:SupBound}
\lim_n \frac{\ln S(n,q)}{\ln n} \as{<} \lim_n \frac{q a_{m+1}(n)-\psi(a_m(n))}{\ln n}.  
\end{equation}

The lower bound is more difficult to obtain and the method proposed here uses the Bienaym\'e-Tchebitchev inequality and is thus only applicable to $\iid$ variables.
 \begin{equation}
 \Prob{ N_m > \frac{\Esp N_m} {2}} < \frac{4 }{n \Esp \chi_{I_m} } 
\end{equation}

 \begin{align*}
 n \Esp \chi_{I_m} &> n \Esp \chi_{I_{R-1}} \\
  &> n \int_{y^\dagger(1-\frac{1}{R})}^{y^\dagger(1-\frac{1}{2R})} e^{-\psi(y)} dy.\\
  &> n \frac{y^\dagger} {2R} e^{-\psi(y^\dagger(1-\frac{1}{2R}) )} \\
  &= \frac{y^\dagger} {2R} (\ln n)^{2\rho(1+\epsilon(y))} \\
  &= (\ln n)^{(2\rho-1-1/\rho)(1+\epsilon(y))} 
\end{align*}
This latest equation imposes the most restrictive condition on $R$. 
Indeed, this implies that
\begin{equation*}
 \sum_r \Prob{ N_m > \frac{\Esp N_m} {2}} < +\infty. 
\end{equation*}
The Borel-Cantelli Lemma can be used to obtain
 \begin{equation}
 \frac{1}{2} e^{q a_{m}-\psi(a_{m+1})} \as{<} S(n,q)  .  
\end{equation}
 \begin{equation}
 \label{eq:InfBound}
\lim_n \frac{q a_{m}-\psi(a_{m+1})}{\ln n} \as{<} \lim_n f(n,q)  .  
\end{equation}

The last needed step is to verify that the two bounds converge when $n\rightarrow +\infty$.
It is needed to verify first that
\[ \frac{q y^\dagger} {R \ln n} = \frac{q y^\dagger \ln \ln n}{(\ln n)^2}. \]
The convergence comes at the price of a condition on $q$
\[ \exists \epsilon>0, q \ll (\ln n)^{2-1/\rho-\epsilon}. \]
The second term is verified easily
\[ \frac{ \psi'(y^\dagger) y^\dagger} {R \ln n} \sim \frac{\ln \ln n}{\ln n}. \]
Thus, the image of $I_m$ by $y \mapsto q y - \psi(y)$ converges towards a simple point. 

Combining these results with (\ref{eq:InfBound}) and (\ref{eq:SupBound}) leads to
\begin{equation}
\label{eq:app:LogConv:Core}
\lim_n \frac{\ln S(n,q)}{\ln n} \as{=} \lim_n \frac{q y^*(n)-\psi(y^*)}{\ln n}, 
\end{equation}
which hence, combined to a saddle point evaluation of the truncated moments proves Eq.~(\ref{eq:SMT}).

\section{Estimation of $\theta$: theoretical results}
\label{sec:theta_appendix}

 \subsection{Construction}
 \label{sec:theta_appendixa}
In order to prove the consistence and the asymptotic normality of $\Omega_k$, it is easier to consider variables $(G_1, \dots , G_n)$ of law $\Lambda^{(k)}$. 
If we then define 
\begin{equation} \Delta_i = i ( G_i - G_{i+1}),\end{equation}
one can see that
\begin{equation} 
\label{eq:app:omega:dist}
p_{\Lambda^{(k)}} ( \delta_1, \dots , \delta_{k-1}, g_k ) = \exp\left(-e^{-g_{k}}-k g_k - \sum_{i<k}i\cdot\Delta_{i}\right).\end{equation}

In other words, $\Delta_i$ and $G_k$ are independent. 
Moreover, if we call 
\begin{equation} \zeta(s;n)= \sum_{j=1}^n \frac{1}{j^s} \end{equation}
Some computations then leads to : 
\begin{equation} \Esp{ \Delta_i}= 1, \end{equation}
\begin{equation} \Var{\Delta_i}= 1.\end{equation}
Similarly with $\gamma$ the Euler-Mascheroni constant:
\begin{equation} \Esp{G_{k}}=\gamma-\zeta(1;k-1),\end{equation}
\begin{equation}\Var{G_{k}}= \frac{\pi^{2}}{6}-\zeta(2;k-1).\end{equation}
 If we now try to apply $\Omega_k$ to the $G_i$ variable, we can translate
\begin{equation} \Omega_k= \sum \alpha_i G_i  \end{equation}
into 
\begin{equation} \Omega_k= G_k + \sum_ {i=1}^{k-1} \beta_i \Delta_i.\end{equation}
with $\beta_i= \frac{1}{i} \sum_{j \leq i} \alpha_j $. We can thus calculate the quadratic standard error  with
\begin{equation} \Esp{\Omega_k^2}= \Esp{\Omega}^2 + \Var{\Omega}. \end{equation}
which brings forth 

 \begin{equation}
\Esp{\Omega_k^{2}}= \frac{\pi^2}{6} -\zeta(2,k-1) + \sum_{i=1}^{k-1} \beta_i^2  +  \left(\gamma -\zeta(1;k-1) + \sum_{1}^{k-1} \beta_i \right)^2 . \end{equation}

We can then minimize $\Esp{\Omega_k^{2}}$ directly according to the $\beta_i$ with a Lagrange multiplier $\lambda$ in order to insure than $\Esp{\Omega_k}=0 $
\begin{equation} \partial_{\beta_i} (\Esp{\Omega_k^{2}}- \lambda{\Esp{\Omega_k}} ) = 2(\gamma-\zeta(1;k-1))-\lambda+2\sum_{j=1}^{k-1}\beta_{j} + 2 \beta_i .\end{equation}
An obvious solution is then $\beta_i= \beta=\frac{\zeta(1;k-1)-\gamma}{k-1} $ 
and so : 
 \begin{equation} 
\begin{cases}
\forall i, i\neq k & \alpha_i= \frac{\gamma- \sum_{l=1}^k \frac{1}{l}}{k-1} \\
 & \alpha_k=1- (k-1)\alpha_1 \\
\end{cases}
 \end{equation}
 
 \subsection{Consistence}
 \label{sec:theta_appendixb}
 We can now verify that $\lim_k \Esp{\Omega_k^2}=0$ for the variables $G_k$
 
 \begin{equation} \Esp{\Omega_k^2}= \frac{\pi^2} {6}  - \zeta(2,k-1) + \frac{(\zeta(1;k-1)-\gamma)^2}{ k-1}. \end{equation}
 The result comes immediately from the fact that $\zeta(2)= \frac{\pi^2} {6}. $
Using the convergence in distributions of the $U_k$, we then obtain

\begin{equation} \Esp{\Omega^2(U)} \underset{k}{\rightarrow} 0 .\end{equation}
 
 \subsection{Asymptotic Normality}
  \label{sec:theta_appendixc}

One can write $\Omega_k$ as 
\begin{equation} \Omega_k= g_k + \beta \sum_ {i=1}^{k-1} \Delta_i.\end{equation}
From Eq.~\ref{eq:app:omega:dist}, it is easy to verify that the $\Delta_i$ variables are $\iid$.
Thus,
\begin{equation} \frac{\Omega_k - g_k -\zeta(1;k-1)+\gamma}{ \zeta(1;k-1)-\gamma } \convD \Normal{0}{1}. \end{equation}
Moreover,
\begin{equation} \Esp {g_k -\zeta(1;k-1)+\gamma}=0, \end{equation}
\begin{equation} \Var {g_k -\zeta(1;k-1)+\gamma} = \frac{\pi^2}{6} - \zeta(2;k-1), \end{equation}
implies that 
\begin{equation} g_k-\zeta(1;k-1)+\gamma \convP 0. \end{equation}
It is thus possible to eliminate this term from the previous equation:
\begin{equation} \frac{\Omega_k}{ \zeta(1;k-1)-\gamma } \convD \Normal{0}{1}.\end{equation}
Combining this result with the convergence in distribution of the $U_k$ show that there is a $k(n)$ such that $\Omega_k$ is asymptotically normal.

\end{appendix}

\end{document}